# The Geometry of Motions

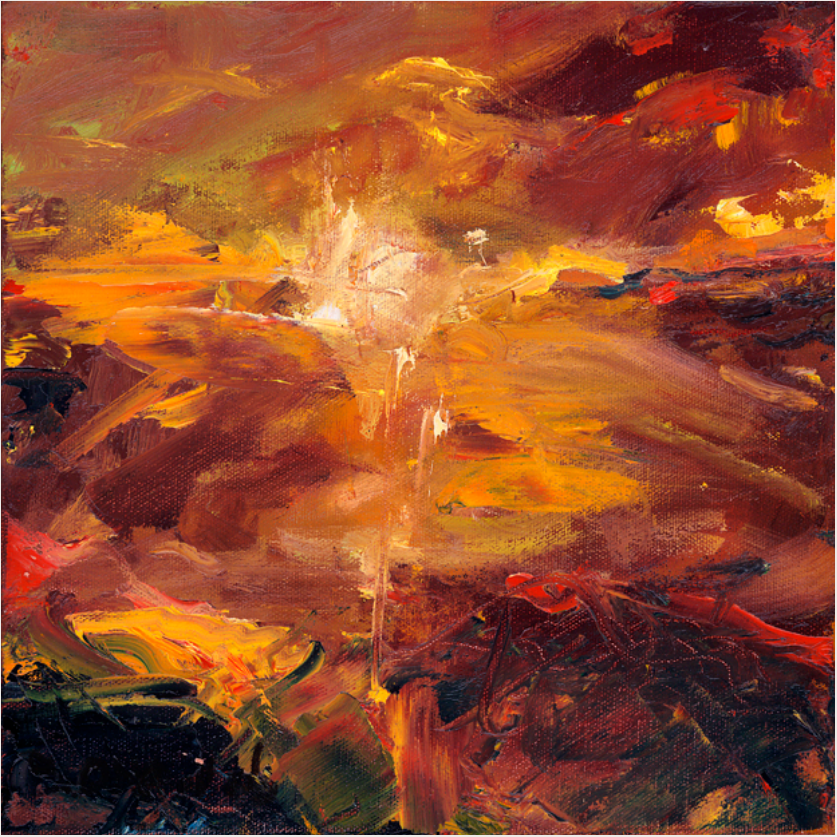

## Patrick Iglesias-Zemmour

# THE GEOMETRY OF MOTION

PATRICK IGLESIAS-ZEMMOUR

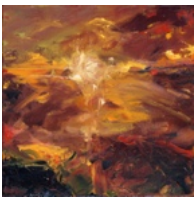

Artwork, *Bereshit.*
Barbara Longwill, Abstract Expressionist.
Oil on canvas 12x12, 2006
`http://www.artlongwill.com`
`artlongwill1@aol.com`

*To six million stars*

# Preface

This work grew from an invitation to speak at the colloquium *La reconquête de la dynamique par la géométrie après Lagrange* (The reconquest of dynamics by geometry after Lagrange), which took place from March 24 to 26, 2010 at the IHÉS in Paris. I cannot begin without first thanking those who were at the origin of this project : Jean-Pierre Bourguignon, Pierre Cartier, and Yvette Kosmann-Schwarzbach. Their invitation provided the initial impetus to formalize these ideas.

The writing of this text was made possible by the kind invitation of Thibault Damour and the hospitality of the IHÉS, whose stimulating environment was invaluable.

I also wish to extend my gratitude to Georges Hansel, for our illuminating discussions on Maimonides' *The Guide for the Perplexed*, which resonated deeply with my own readings. Finally, my thanks to Jean-Jacques Szcecyniarz for his insightful remarks and his encouragement to publish this work.

This volume, and the trilogy it begins, is the result of these conversations and the generous support of these individuals and institutions.

The author acknowledges the significant help of the AI assistant Gemini (Google) in the collaborative refinement and editing of this manuscript, without which this book would not have been possible in its present form.

# Introduction

> « *Since nature is a principle of motion and change, and our inquiry is concerning nature, what motion is must not escape us; for, if it is unknown, nature must necessarily be unknown.* »

> — Aristotle, *Physics*, Book III

Aristotle's ancient charge remains the fundamental task of physics. This work takes that charge as its central axiom : to understand nature is to understand motion. It proposes that a physical theory of motion is, in its essence, a geometry in the sense of Felix Klein, from which the dynamical laws themselves emerge as consequences. The fundamental character of this geometry is defined by its **inertia group** —the group of spacetime transformations that preserves its class of most elementary, "in-articulated" motions.

This first volume embarks on a journey through the three great **epistemological ruptures** that have defined our understanding of motion. We will trace the dissolution of successive worldviews by translating the foundational principles of physics into the unambiguous language of group theory. Vague philosophical statements such as "space is absolute" or "time is relative" will be recast as precise, falsifiable mathematical propositions about the structure of spacetime and the symmetries it admits.

Our path begins with the mechanics of Aristotle, where the seemingly intuitive categories of absolute Space and absolute Time are formalized. We will construct the **Group of Aristotle** as the unique group of symmetries preserving this rigid structure. In this framework, the inertial motions are simply states of rest. Far from being merely a historical artifact, we will see that this Aristotelian geometry provides the essential baseline against which the dissipative forces of our sublunar world can be measured.

We will then witness the collapse of this edifice, not through a single discovery, but through a series of profound objections that form the twin cornerstones of modern mechanics. We will examine the incisive philosophical critique of Maimonides on the nature of time as an "accident of motion," a critique that consciously sided with Plato against Aristotle's concept of an eternal, disembodied Time. In parallel, we will explore the revolutionary thought experiments of Giordano Bruno and Galileo, which demonstrated the physical indiscernibility of rest and uniform rectilinear motion. This second rupture, driven by the abandonment of Aristotelian Space, leads to the construction of a new, more subtle geometry defined by the **Galilean Group**.

The final rupture is precipitated by the crisis in nineteenth-century electromagnetism. The work of Lorentz, Poincaré, and Einstein reveals that even the Galilean structure is untenable. The constancy of the speed of light forces a fusion of space and time, finally vindicating Maimonides' ancient objection. This leads to the geometry of the **Poincaré Group**, where both absolute Space and absolute Time are lost, replaced by the invariant structure of the Minkowski quadratic form.

**About This Volume.**   This book provides a self-contained, mathematical account of these historical transformations. By focusing on the underlying inertia groups, it offers a clear and rigorous framework for understanding the transition from the classical to the modern conception of spacetime. The central theorems of this volume demonstrate with mathematical certainty why there can be no Galilean "Space" and no Einsteinian "Time" in the absolute sense that Aristotle conceived them. It is an argument for seeing mechanics not as a collection of disparate models, but as a unified story of evolving geometric principles, each built upon a different, logically consistent axiom of inertia.

**A Foundation for a Larger Project.**   Yet, this historical deconstruction is only the first step. The abandonment of absolute Space and Time forces us to confront motion itself as the primary object of study. The true stage for physics is not spacetime, but the infinite-dimensional **space of all possible motions**.

This volume serves as the indispensable philosophical and mathematical prologue to a larger project. The second part will develop the theory of **dynamics**, understood as the art of comparing real motions to ideal ones. We will see how this principle, initiated by Newton and perfected by Lagrange through his method of "variation of constants," reveals that this comparison is governed by a canonical **symplectic structure** on the space of motions. Finally, the third volume will show how the newcomers to physics —quantization and field theory— demand a corresponding evolution in our geometric tools. To describe the infinite-dimensional or singular spaces that these theories inhabit, the very notion of a smooth space must be enlarged, completing the journey from the classical world to the frontiers of physics. This book, therefore, is the first part of a trilogy dedicated to "The Geometry of Motion," establishing the classical foundations and the group-theoretic language necessary to embark on the study of dynamics in its most fundamental form.

# THE GEOMETRY OF MOTION

## MECHANICS AS GEOMETRIES

« *Therefore take this portrait as it is, for it is not addressed to you to instruct you in what you already know, nor to add a little water to the torrent of your intelligence and judgment; but because it is not our custom, as far as I know, to disdain the portrait and presentation of things, even if we know them better from life.* »

Giordano Bruno, *The Ash Wednesday Supper* (Le Banquet des Cendres).

# 1. The Principle of Motion

The very principle of nature is motion, this is what Aristotle asserts in *Physics*[1]. Consequently, any reflection on nature is first a reflection on motion. This is how the first chapter of book III begins :

> « Since **nature is a principle of motion** *and change, and our inquiry is concerning nature,* **what motion is must not escape us***, for, if it is unknown, nature must necessarily be unknown.* »

To apprehend motion, describe it, and then determine what follows from it, Aristotle deconstructs it, and this leads him to abstract the category of *Time* and that of *Place* (τόπος), which for the purposes of physics is synonymous with the modern idea of *Space*. Indeed, it is through these *categories* that the motion of things is ordinarily expressed.

This abstraction of Time and Space into absolute, pre-existing categories was Aristotle's foundational choice. It was not, however, the only view. Plato, his teacher, had argued for a different conception, one in which time was not an eternal container for motion, but an emergent property of the cosmos, created with the heavens themselves. [2] This ancient debate between two competing axioms would lie dormant for centuries, but its resolution is central to the story of modern mechanics.

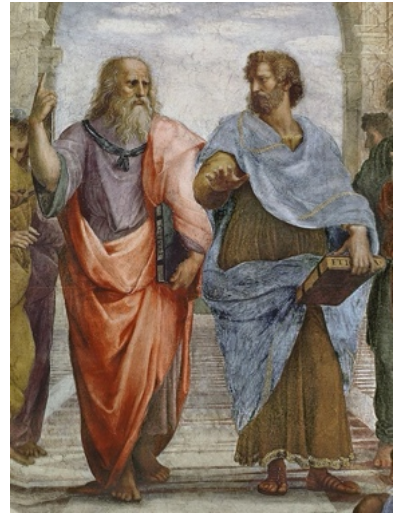

*Plato* and *Aristotle* at the Lyceum. *Their foundational disagreement on the nature of Time —whether an absolute container for motion or an accident of the cosmos— established a central debate that shaped the next two millennia of physics.*

Thus, during the centuries that followed, Aristotle's Time and Space were considered by physicists as the two *imperative* categories of discourse on nature, without which no motion is imaginable or describable. They conditioned and shaped thought to the point that they still weigh on the minds

---

of physicists, long after the Galilean and relativistic revolutions, which nevertheless refuted these postulates.

But does this mean that there is neither time nor space? No, obviously not, but the time as we feel it, whether it passes quickly or slowly depending on what we are experiencing, or the space as we conceive it, delimited by the walls surrounding us and the roof above us, are not the Time and Space that Aristotle desires, as he defines them.

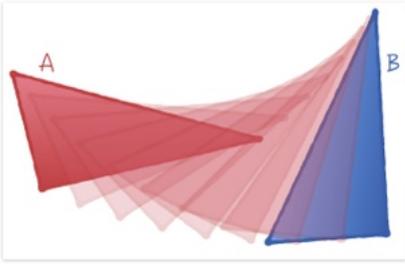

**Transformation Groups** The case of triangle congruence is expressed by the existence of a transformation, composed of a rotation and a translation, that maps triangle A onto triangle B. These transformations *group together* into what is called a *group* of transformations. Mathematically, groups are sets of transformations that can be composed with each other and also inverted.

To translate these foundational principles into a precise and falsifiable form, we will employ a powerful mathematical language. We propose that the essence of each mechanics —its core set of axioms about the nature of motion— can be perfectly captured by its **Inertia Group**. This is the group of all spacetime transformations that preserve the system's privileged class of inertial motions. This Inertia Group becomes the unique signature of the mechanics, embodying its fundamental symmetries and defines its **Geometry**.

It is this passage to the concept of the *Inertia Group* that clarifies the exposition of the different epistemological breaks, by mathematically distinguishing them, and which leads us, over the centuries, and millennia, from Aristotelian mechanics to Einsteinian mechanics. It is this group that completely encapsulates the internal structure of the mechanics considered, eliminating any ambiguity, and rendering vain any allusive discourse about their properties.

## 2.  Places or Aristotle's Space

Aristotle dedicates chapter 4 of book IV of *Physics* to the category of *place* (τόπος), which, along with time, is one of the imperative categories of physics. Here is how he introduces this chapter :

> « *One must first understand that one would not research about place if there were no **motion according to place**.* »

It is obviously because motion according to place is trivially described by the list of successively occupied places (and by the list of durations during which each place is occupied) that it is important for Aristotle to define the notion of place, to specify its characteristics and delineate its aspects. Here is what he says about it, in the continuation of the chapter. Firstly :

« *Place is the **first containing** of that of which it is the place.* »

Which can be read as follows : there is nothing prior to the place in terms of containing things, *i.e.*, there is nothing between the place and the thing other than the place and the thing. The thing moves continuously [3] from place to place. Secondly :

« *It is **nothing of the thing**.* »

The place is not confounded with the thing, *i.e.*, the place exists *per se* absolutely, the place is objective, it is not an *accident* of the thing. Furthermore, it is admitted that any place in this world is susceptible to being, at one time or another, the place of rest for some thing. (One need only imagine an object constructed such that its natural equilibrium, its "proper place," is precisely there.)

The set of all possible places of rest for all things, taken as a whole, are what physicists subsequently called *space* [4].

Thirdly and finally, once what place is or is not has been clarified, Aristotle then adds :

« *We certainly consider that* […] *each of the bodies is carried by nature and remains in its **proper places**.* »

We interpret the proper places of a thing as its *places of rest*, the places where its potential motions cease : a stone that rolls and stops in the hollow of a path, a drop of dew on the edge of a leaf, etc. In Aristotelian mechanics, everything tends towards rest. This is Aristotle's axiom of ***absolute rest***. By admitting that every place can be the place of rest for a thing, we can represent **Aristotle's Space** as the set of places of rest for things, or equivalently as the set of *resting motions*. Thus :

> *The central axiom of Aristotelian geometry is established :*
> ***Aristotle's Space is the set of all resting motions.*** *This identification is the crucial step that allows us to analyze this mechanics within the language of motion itself.*

---

3. Indeed, in the same paragraph, one reads « *Motion seems to belong to the continua.* »

4. It is customary for a structured set to be called a *space*, in this case, the set of places according to Aristotle has a Euclidean structure : we know how to measure distances between places and angles between directions to these places.

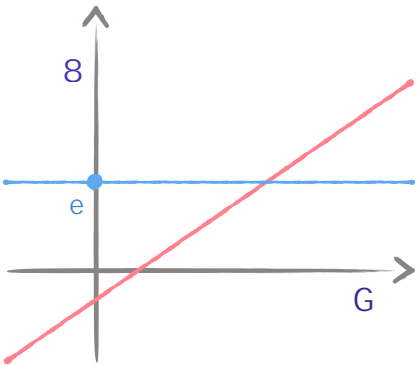

**Resting Motion and Uniform Motion.** *The blue line is a rest at place r, and the red line is a uniform rectilinear motion.*

**Aristotle's Rests.** The expression *resting motion* resembles an oxymoron [5], and yet... An Aristotelian motion is a map from Time to Space. Rest appears in this way as a horizontal line of height *r* in the product Space × Time, where *r* denotes the place of rest. On the other hand, uniform rectilinear motion is described by an inclined line. Rest can then be seen as a uniform rectilinear motion whose peculiarity is having zero inclination. Rest is thus a motion, in the same way that zero is a number. This should not surprise us ; we have long been accustomed, in science and in mathematics in particular, to representing *the absence of something as the presence of nothing*.

**Aristotle's Inertia.** It is this identification, between place, place of rest, and resting motion, that will allow us later to explain the *group of Aristotle*, *i.e.*, the **Inertia Group** of Aristotelian mechanics. We will see, in the following chapters, how to each different mechanics (Aristotelian mechanics, Galilean mechanics, and Einsteinian mechanics), is associated, on the other hand, a set of distinguished motions, which we call ***inertial motions***, and on the other hand, a certain group of spacetime transformations, the Inertia Group, which preserves, among other properties, the set of inertial motions. The inertial motions of Aristotelian mechanics are, as they should be, the resting motions. The Inertia Group captures the essence of the mechanics it represents. The word inertial is to be taken in its etymological sense *in-art → inert*, *i.e.*, an *in-articulated* motion, without spring, minimal, passive...

---

5. It would be just as legitimate to simply call these motions *rests*.

Contrary to what one might imagine, the use of the adjective inertial in the Aristotelian framework is entirely legitimate, even if it does not coincide with that ordinarily used in Galilean mechanics.

> INERTE, lat. *iners*, sans activité, sans ressort ; il y a l'inertie du corps et celle de l'esprit ; force d'inertie, résistance passive ; le lat. *iners*, composé de *in* nég. et *artus*, membre ; inerte, qui n'est plus articulé en quelque sorte.

**Etymological Dictionary
of the French Language**
Adolphe Mazure, 1863.

Indeed, a certain laziness of the mind, born from the revolutionary surprise of the Galilean principle, has tended to obscure a crucial point : the state of rest played the exact same conceptual role for Aristotle that uniform rectilinear motion played for Galileo. It was the fundamental, unforced state—the very definition of inertia in his system. Inertia does indeed exist in Aristotle ; it is rest, immobility, a concept entirely consistent with his vision of a central, immobile Earth.

We will see in the second volume of this essay, how the existence, in all cases, of a family of distinguished motions — playing the role of inertial motions, although they may no longer necessarily be rectilinear or uniform — was the basis for Joseph-Louis Lagrange's first symplectic constructions. He pioneered this approach by treating the real, perturbed motion of a planet as a curve within the space of ideal Keplerian orbits, effectively promoting the ellipses themselves to the status of ideal motions. This transcended Space and Time, towards the category of Motion.

**The etymology of the word *inert*** and its derivatives : *inertia, inertial*, is interesting for more than one reason. As we noted above, the word *inert* is a deformation of the composition *in-art*, where the prefix *in* is privative. The word inert therefore means *without art.* The word *art* itself has a distant origin and significant ramifications.

According to the attached excerpt from Jauffret's etymological dictionary, the root of this word is the mimologism *Re, Roo, Aroo*, which describes the action of ploughing—the first art of men—expressed in Old French by the verb ***arer*** (to plough) and the noun ***arage*** (ploughing). From these roots, we derive the Old French words ***are*** (still in use in modern French with ***hectare = hundred ares***) a unit of land area (akin to the English *acre*), and *aire*, which becomes *area* in contemporary English, and many others in the same vein.

ART, profession, occupation quelconque. Ce mot, qui a un si grand nombre de dérivés, appartient à une famille extrêmement nombreuse. Il faudroit une dissertation entière pour développer son origine et ses différentes ramifications. Il remonte au mimologisme *re*, *roo*, *aroo*, qui peint l'action de labourer. L'*arage*, ou l'action de faire aller la charrue, a été le premier des arts, l'art par excellence. Dans les langues orientales, *ar*, *arz*, signifie la *terre*. De-là aussi le latin *arvum*, un champ, *area*, plaine, campagne. ( *Voy.* TERRE.) *Aroo* a produit le mot latin *armus*, le bras, que les Allemands et plusieurs autres peuples du Nord ont conservé. *Armus* a produit *arma*, les armes; *armilla*, bracelet, cercle; d'où nous avons fait sphère *armillaire*, composé des cercles astronomiques.

**Etymological Dictionary**
**of the French Language**
Louis-François Jauffret. Paris, An VII.

This root is also related to oriental languages, by the word *arz*, in Hebrew for example ארץ, consisting of the letters *a.r.tz*, and pronounced *äretz*. This same root gives the Latin word *armus* from which *arm* is derived, since one ploughs by the strength of one's arms. We also have *army*, *armilla* (a bracelet), etc. From this proximity to the arm, we derive *articulated, article, articulation, the French orteil* (toe), etc. We can thus see that the use we make of the word *inertial* conforms to its origin in all cases : Inertial motions are unarticulated; they are not the fruit of any particular art.

## 3. Aristotle's Time

*Following what has just been said, we must proceed concerning time* (χρόνος). It is thus, in chapter 10 of book IV, that Aristotle introduces his reflection that will lead him to establish an absolute, eternal, disembodied Time. We read, in chapter 11 of book IV of *Physics*, the elements and principles that structure the Aristotelian conception of time, and which will make it the second imperative category of physics regarding motion :

« *This is what time is : the **number of motion**, according to earlier and later.* […] *Time is the counted, not that by which we count.* »

We also read in chapter 12 :

« *Time is the **measure of motion**.* »

And to measure the motion of things, nothing is better than comparing it to the perfect motion of the spheres, as said in chapter 14 :

> « *If then that which is primary is the **measure of all things** of the same kind, **uniform circular translation is primarily the measure** because its number is the most knowable.* »

Time is thus abstracted from uniform circular motion, precisely :

> « *This is why time **seems to be the motion of the sphere**, because **by this motion also other motions are measured** and also time.* »

Thus, Aristotle begins by linking time to motion, a principle consistent with Plato's view that time was created with the heavens and is inseparable from their movement. But he then makes a profound and decisive axiomatic leap. By postulating the existence of a perfect, eternal motion —that of the celestial sphere— he proposes that its time can be detached from it to become a universal standard. The chronology of the sphere becomes **Aristotle's Time**, an absolute and disembodied entity to which the motion of every other thing must be referred. This numbering of the sphere's motion (the *quantity* of motion of the sphere) is precisely what sundials offer us. It should also be noted that, since Aristotle admits that the motion of the sphere is eternal, Aristotle's Time is both absolute and eternal.

The critique of simultaneity —that is, if the fact that two distant events occur "at the same time" is consistent and physically verifiable— is not considered, and it cannot be, as such a question would have been absurd in the context of that ancient time. But we will see that it is this question of simultaneity that will later invalidate Aristotle's principles on time, and we will see the consequences of such an abandonment later. For now, the existence of a *perfect eternal motion* to which all other motion can be related, and the abstraction of this eternal motion from the thing whose motion it is, is what is meant by Aristotelian Time. This requires much more from nature than Plato's time, which was created with the "Heaven." [6].

This concept of a universal, absolute Time provides the second great pillar of the Aristotelian world. Just as his Space is a universal stage of absolute places, his Time is a universal river flowing at the same rate for all things. Together,

---

6. Plato, *Timaeus*. Translated by Donald J. Zeyl. Hackett Publishing Company, 2000. « However it be, time came into being together with the heaven » they form an ideal [7] rigid, absolute framework within which all motion is to be measured. We have now established the two foundational principles of his mechanics. The next step is to translate this intuitive worldview into a precise mathematical structure : its Inertia Group.

## 4. The Group of Aristotle

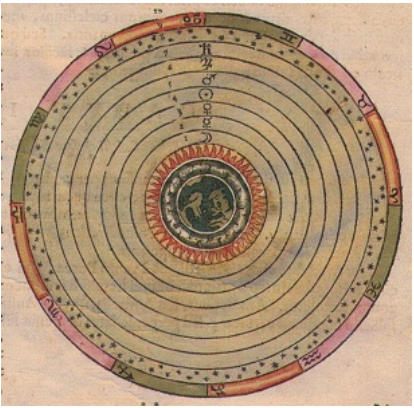

**Aristotle's Cosmos**
*The Earth is fixed at the center of the universe, the celestial sphere gives the measure of time.*

Aristotle's Space and Time have a mathematical translation that is no longer a mystery today. They are represented by *oriented Euclidean affine spaces*, which we will denote E and T respectively. Space E is of dimension 3 : once an origin of space and three base points are chosen, one in each direction, three numbers are needed, the *coordinates*, to unambiguously locate any point in E. Space T is of dimension 1 : once the origin of time and a base point are chosen, a single number is sufficient, *the date*, to locate an instant in the continuum of Aristotelian Time.

Distance in E and T is defined by the Pythagorean theorem, which states that the distance of a point from the origin is the square root of the sum of the squares of the coordinates of that point. The standard model for E is the space $\mathbf{R}^3$ of triplets of numbers $r = (x, y, z)$. The model for T is simply $\mathbf{R}$. The set of pairs $q = (r, t)$ describes the product space E × T, which is generally called ***Spacetime***.

It can already be noted that the Euclidean structure of T is naturally inherited from the Euclidean structure of E, insofar as time is read on sundials, or dial clocks, or clepsydras, or by any other mechanism that traces the motion of the sphere by the position of hands or markers.

We have previously stated that we will identify the places of things with their proper places, *i.e.*, with their places of rest, and then the places of rest with resting motions. How is this then translated in the Spacetime scheme we have just described?

---

7. Some say that the (mathematical) straight line is *what you have in mind when you draw a chalk line on the blackboard.*

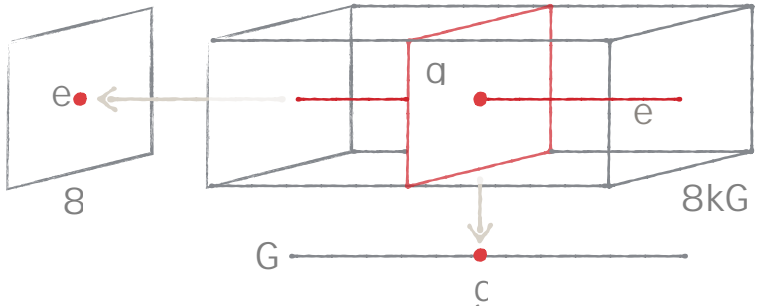

Aristotle's Spacetime

**Statement.** *Space, the set of places according to Aristotle, is identified with the set of resting motions of* material points*, i.e., with the set of graphs of constant maps* $[t \mapsto r]$*, from* $\mathrm{T}$ *to* $\mathrm{E} \times \mathrm{T}$*, where* $r$ *spans* $\mathrm{E}$*. They will be denoted by a bold letter*

$$\mathbf{E} = \big\{ \mathbf{r} = \{(r, t) \mid t \in \mathrm{T}\} \subset \mathrm{E} \times \mathrm{T} \;\big|\; r \in \mathrm{E} \big\}.$$

*Elements* $\mathbf{r}$ *are represented by the horizontal red lines in the figure above. Similarly, we identify Aristotle's Time with the sheets of simultaneous events*[8] *, i.e., the sets* $\mathbf{t}$ *in* $\mathrm{E} \times \mathrm{T}$ *of events occurring at the same instant* $t$*, as* $t$ *spans* $\mathrm{T}$

$$\mathbf{T} = \big\{ \mathbf{t} = \{(r, t) \mid r \in \mathrm{E}\} \subset \mathrm{E} \times \mathrm{T} \;\big|\; t \in \mathrm{T} \big\}.$$

*They are represented by the vertical red slices of* $\mathrm{E} \times \mathrm{T}$*, above* $t$*, in the figure above.*

It is clear that the spaces $\mathbf{E}$ and $\mathbf{T}$ preserve their Euclidean structures, inherited from $\mathrm{E}$ and $\mathrm{T}$. We will denote this Euclidean distance by $d$. By definition, if we set $r = (x, y, z)$ and $r' = (x', y', z')$, then

$$d(\mathbf{r}, \mathbf{r}') = \sqrt{(x' - x)^2 + (y' - y)^2 + (z' - z)^2}. \tag{$\clubsuit$}$$

Similarly, concerning time :

$$d(\mathbf{t}, \mathbf{t}') = \sqrt{(t' - t)^2} = |\, t' - t \,|. \tag{$\spadesuit$}$$

---

8. It is customary to call elements of spacetime *events.*

It is the group of Aristotle [9] which we will say is the **inertia group** of Aristotelian mechanics, that structures Aristotle's Spacetime. For, as we will see later, Galilean or relativistic mechanics share the same spacetime as a real affine space. It is therefore important to understand that it is the structure with which this spacetime is endowed by the inertia group that makes the difference.

Precisely, the **Group of Aristotle** is defined as the group of transformations of the spacetime $E \times T$, which preserves the categories of Time and Space, their measure or distance (and possibly their orientation). More precisely, an *Aristotelian transformation* is an automorphism $g$ of $E \times T$, that satisfies the following conditions :

- ✔ **Preserves Space** : $g$ exchanges resting motions. *i.e.*, the image of a resting motion $\mathbf{r} = \{(r, t) \mid t \in T\}$ by $g$ is another resting motion $\mathbf{r}' = g(\mathbf{r})$.

- ✔ **Preserves Time** : $g$ exchanges temporal slices. The image of a sheet of simultaneous events $\mathbf{t} = \{(r, t) \mid r \in E\}$, is another sheet of simultaneous events $\mathbf{t}' = \{(r, t') \mid r \in E\}$.

- ✔ **Preserves the Euclidean structures** of Space and Time : $d(\mathbf{r}, \mathbf{r}') = d(g(\mathbf{r}), g(\mathbf{r}'))$ and $d(\mathbf{t}, \mathbf{t}') = d(g(\mathbf{t}), g(\mathbf{t}'))$.

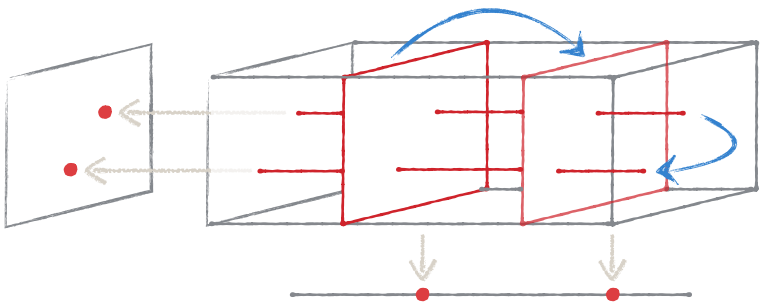

The Group of Aristotle

**Proposition.** *Acting on the spacetime* E × T, *the group of Aristotle is represented by the matrices*

$$g = \begin{pmatrix} A & 0 & C \\ 0 & 1 & e \\ 0 & 0 & 1 \end{pmatrix} : \begin{pmatrix} r \\ t \\ 1 \end{pmatrix} \mapsto \begin{pmatrix} Ar + C \\ \pm t + e \\ 1 \end{pmatrix} \quad where \quad \begin{cases} A \in O(3) \\ C \in \mathbf{R}^3 \\ e \in \mathbf{R} \end{cases} .$$

*Proof.* By definition $g$ acts on spacetime, we can write $g(r,t) = (\rho(r,t), \tau(r,t))$ where $\rho$ and $\tau$ are certain functions. Now, since for every $r$ there exists $r'$ such that $g(\{(r,t) \mid t \in T\}) = \{(r',t) \mid t \in T\}$, it must be that $\rho(r,t) = r'$ for any $t$.

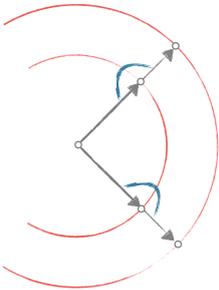

Furthermore, since for every $t$ there exists $t'$ such that $g(\{(r,t) \mid r \in E\}) = \{(r,t') \mid r \in E\}$, then $\tau(r,t) = t'$ for any $r$. Thus $g(r,t) = (\rho(r), \tau(t))$ where $\rho$ is a transformation of E and $\tau$ a transformation of T. For simplicity, we identify E with $\mathbf{R}^3$. Let's now establish that Euclidean transformations of space E can be written as $r \mapsto Ar + C$, where $A \in O(3)$ and $C \in \mathbf{R}^3$.

This is a standard proof but we provide a version to show how these things are articulated. Consider three points $a$, $b$, and $c$ and their images $a'$, $b'$, and $c'$ under $\rho$. By hypothesis we have the three distance equalities : $d(a,b) = d(a',b')$, $d(a,c) = d(a',c')$ and $d(b,c) = d(b',c')$. This is a case of triangle congruence. The angles at points $a$, $b$, and $c$ are therefore equal to the angles at points $a'$, $b'$, and $c'$. The transformation $\rho$ therefore preserves angles. Let $A(r) = \rho(r) - C$ with $C = \rho(0)$.

The map A then preserves the origin and obviously also distances. Let $r$ be a point in E and $\lambda$ a positive number greater than 1. The image $A(r)$ is on the circle of radius $\|r\|$, and $A(\lambda r)$ on the circle of radius $\|\lambda r\|$. Furthermore, the triangle $(O, r, \lambda r)$ with vertex $r$ is flat, since A preserves angles, the triangle $(O, A(r), A(\lambda r))$ with vertex $A(r)$ is also flat, from which we deduce that $A(\lambda r) = \lambda A(r)$. Other cases for $\lambda$ are handled just as simply, and this equality holds for all $\lambda \in \mathbf{R}$.

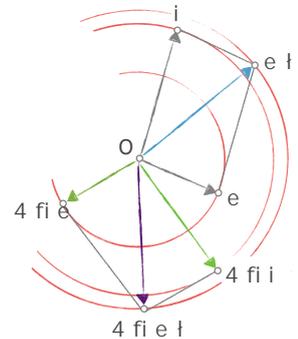

Next, consider two vectors $r$ and $v$, and the triangle with vertices $0$, $r$, and $r + v$. Knowing that A maps every triangle to a congruent triangle (in the sense of triangle congruence cases), the triangle $(r, v, r + v)$ is congruent to the triangle $(\mathrm{A}(r), \mathrm{A}(v), \mathrm{A}(r + v))$, and thus, as the figure shows, $\mathrm{A}(r + v) = \mathrm{A}(r) + \mathrm{A}(v)$. The map A is therefore linear, it preserves distances and angles, it is thus, by definition, an orthogonal transformation. Thus, $\varrho(r) = \mathrm{A}r + \mathrm{C}$, with $\mathrm{A} \in \mathrm{O}(3)$ and $\mathrm{C} \in \mathbf{R}^3$.

Finally, what is said for E is *a fortiori* true for T, and we obtain $\tau(t) = \pm t + e$, with $e \in \mathbf{R}$. $\qquad\qquad\Box$

The action of Aristotle's group is therefore decoupled : on Space it is identified with the group of Euclidean transformations $r \mapsto \mathrm{A}r + \mathrm{C}$, and similarly on Time where we find again the Euclidean transformations $t \mapsto \pm t + e$. Usage dictates that only the identity component is considered. Aristotle's group is therefore given, strictly speaking, by $r \mapsto \mathrm{A}r + \mathrm{C}$ with $\mathrm{A} \in \mathrm{SO}(3)$ and $t \mapsto t + e$.

**Conclusion.** We have not only equipped Aristotle's Space and Time with their Euclidean structures, which is undoubtedly quite standard and certainly not sufficient to capture the principles of his mechanics. It is by translating these principles through the conservation of both resting motions and temporal sheets of simultaneous events, as well as measures of length and duration, that we have constructed the group that faithfully captures these two principles.

With the introduction of this group of Aristotle, mechanics is then interpreted as a geometry in the sense of Klein. Thus, Spacetime for example — the primary object of this mechanics — becomes subordinate to the group of Aristotle, of which it is a homogeneous space. But more than that, and this is what expresses in a precise formal way Aristotle's principles on the absolute nature of space and time :

- ✔ The space **E** of resting motions is a homogeneous space of Aristotle's group.
- ✔ The temporal line **T** of sheets of simultaneous events is a homogeneous space of Aristotle's group.

We can already see how the use of group theory in this context removes the ambiguities of ordinary discourse. Indeed, a vague sentence such as : *For Aristotle, space and time are absolute*, which is often found as is in many textbooks, is formally translated by the proposition : *Space and time, as conceived by Aristotle, are homogeneous spaces of the group of his mechanics*, and no longer suffers from any ambiguity.

It then appears clearly that the geometry, in Klein's sense, defined by Aristotelian mechanics, is ordinary Euclidean geometry. *i.e.*, the geometry of transformations of space **E** and time **T** that preserve lengths and angles, as regards Space, and durations, as regards Time.

**Note.** In this paragraph, we have identified two particular partitions of spacetime whose role proved essential in the construction of Aristotle's group. These are the set of resting motions on the one hand and the sheets of simultaneous events on the other. We will see later that while the existence of an invariant partition into sheets of simultaneous events is invalidated in Einsteinian relativistic mechanics, the identification of an invariant class of motions in each mechanics, Aristotelian, Galilean or relativistic, remains fundamental in the construction of their inertia group. These motions, we call them **inertial motions**, and we naturally have for the present case :

**Definition.** *Aristotle's rests are the* inertial motions *of Aristotelian mechanics.*

## 5. Maimonides' Objection

*The Guide for the Perplexed*[10] is a major work by Maimonides. He wrote it to reconcile 12th-century CE Jewish intellectuals with their history, their traditional culture, their rites, and their faith. Written for his student Rabban Joseph ben Judah, it is actually addressed to all Jewish scholars whose knowledge and learning — reason — is sometimes challenged by biblical narratives, when their terms are taken literally.

---

10. Maimonides, *The Guide for the Perplexed*. Translated by Shlomo Pines. University of Chicago Press, 1963.

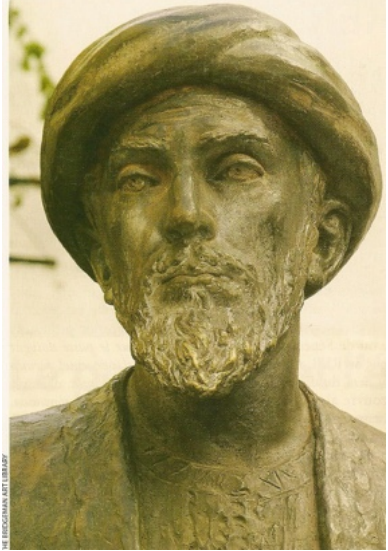

**Maimonides**,
Cordoba 1138 — Fostat 1205.

Maimonides is an admirer of Aristotle, whom he himself calls *the prince of philosophers*. He generally adopts his philosophy. His approach is not one of dogmatic rejection, but of rigorous rational analysis, demanding that scientific claims—even from Aristotle—be supported by irrefutable demonstration. It is therefore not a question for him of disputing Aristotle's heritage, of which he also feels a holder, but neither is it a question of abandoning his own cultural heritage. The strategy of reconciliation between tradition and science that Maimonides proposes, consists simply in seeking the hidden meaning — allegorical — of the sentences and parts of the discourse in the Torah [11], or other sacred texts, when they contradict the reason and knowledge of the time. This leads Maimonides to detailed explanations on the nature of anthropomorphic terms when attributed to the divine, and to give a reading code for sacred Hebrew texts that aligns with Aristotelian concepts and scientific reality. But we will see that when it comes to time, Maimonides can only distance himself from Aristotle. We read this in the introduction to Part II of the *Guide* :

> « *To the preceding propositions I will add one* [...] *which **Aristotle claims to be true** and most admissible ; we concede it to him as a hypothesis, until we*

---

11. The Torah, in Hebrew תורה, consists of the five books of Moses. The word is generally translated as *The Law*. It derives from the word *tor* (תור) which denotes a sequence, a succession ; in this case, a succession of commandments. It is the first part of the Tanakh, in Hebrew תנ"ך — the Hebrew Bible — which consists of : The Law, the Writings, and the Prophets.

> *have been able to set forth our views on this matter. This proposition* [...] *says that **time and motion are eternal** and always existing in act. From this proposition, it necessarily follows that **there is a body that has eternal motion**, always in act, and this is the fifth body*[12]. »

Indeed, this is what Aristotle affirms in chapter 1 of book VIII of *Physics* : « *It is clear that motion is eternal* ». He knows perfectly well that this opinion is subject to disagreement, since he immediately adds : « *Plato alone conceives of time as generated : he says that it was born with the heaven, and that the heaven was born* ».

As can be seen, the stake for Maimonides here is significant : it is the very principle of *Creation*, a cornerstone of Jewish tradition, as shown by the first words of the book of creation in the sacred scriptures. But he also expresses himself clearly, in the *Guide*, chapter XXV Part II :

> « *Know that, **if we avoid professing the eternity of the world, it is not because the text of the Torah would proclaim the world created***; for the texts indicating the creation of the world are no more numerous than those indicating God's corporeality. Regarding the creation of the world too, **the means of an allegorical interpretation would not fail us** nor be forbidden to us; we could use this mode of interpretation here, as we did to reject corporeality (of God)* [...] *But two reasons have led us not to do so. The first is this : God's incorporeality has been demonstrated, and it is necessary to resort to allegorical interpretation whenever the literal meaning is refuted by a demonstration. **But the eternity of the world has not been demonstrated**, and, consequently, **it is not appropriate to do violence to the texts and***

---

12. The *fifth essence*, Aristotle's clock.

> **interpret them allegorically, to uphold an opinion**
> **whose opposite could also be upheld, by reasoning**
> **in another way**. »

Maimonides' position, therefore, is the opposite of dogmatism. He explicitly states that his adherence to the literal text is conditional on the lack of a compelling scientific demonstration to the contrary. He simply considers that the arguments of Aristotle, or his successors and epigones, on the eternity of the world are not a demonstration; and that it is therefore not necessary to renounce the literal interpretation of the texts. Maimonides has already accustomed us to his minimalist philosophical positions. All this has significant repercussions on his conception of time and his disagreement with Aristotle on the subject. Which he explains in chapter XXV of the second part of the *Guide* :

> « *My aim, in this chapter, is to set forth that* **Aristotle**
> **has no demonstration for the eternity of the world**
> *according to his opinion.* **He is not even mistaken**
> **about it**; *I mean, he himself knows that he has no*
> *demonstration for it.* »

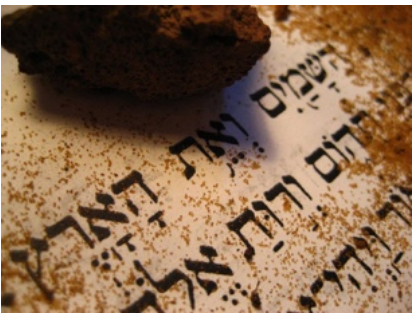

בראשית ברא אלהים את השמים ואת הארץ
*In the beginning Elohim created*
*the heavens and the earth.*

Moreover, knowing the requirements of a *demonstration*, Maimonides does not claim to want to demonstrate that the world was created (as seen above). He merely challenges what is attributed to Aristotle as a demonstration of the eternity of the motion of the celestial sphere; and therefore of the world. Maimonides' argument rests essentially on the simple idea that what has reached its final stage of evolution does not exclude that it was in a different previous state, even if its final stage seems *perfect* and *imperishable*. In other words, if we are unable to imagine a beginning and an end to the motion of the celestial sphere, this does not imply that it was not created. As this debate is only parallel to the question we are addressing here, we will limit ourselves to quoting this passage from chapter XVII, second part of the *Guide*, which deals with that :

> « *Aristotle sets himself against us, arguing against us from the nature of being arrived at its final, perfect state and existing in act, while we,* **we affirm to him that after having arrived at its final state and become perfect, it resembles in nothing what it was at the moment of coming into being**, *and that it was produced from absolute nothingness.* »

It is precisely this intense intellectual dispute with Aristotle that forces Maimonides to go further than Plato and to formulate with new conviction and precision what time truly is. He concludes that its nature is so obscure, even to scholars, for a specific reason : it is not a fundamental substance. It is merely an "accident of motion," which is itself only an "accident of the moved thing." We read in *The Guide*, Part II, chapter XIII, his remarkable explanation :

> « *What we say is that* **time** *has remained an obscure matter for most men of science, so that they have been undecided* [...] *on the question of whether it has, or does not have, real existence, is that it* **is an accident** *in another accident.* **But accidents whose substrates are other accidents** [...] **are a very obscure matter**, *especially when it is joined to this circumstance that the accident which serves as a substrate is not in a fixed state, but changes condition ; for then the matter is more obscure. Now in time, both things are combined ; for firstly* **it is an accident inherent in motion, which is an accident in the thing moved**. »

Maimonides puts time back in its place : *it is only an accident of motion, which itself is only an accident of the thing*. He can only adhere to Plato's view, which is also his own and which is Aristotle's starting point before he elaborates, as he recalls in the introduction to Part II of the *Guide* (fifteenth proposition) :

> « *Motion exists only in time, and* **time cannot be conceived except with motion**. »

Rejecting the eternity of the world proposed by Aristotle, Maimonides can only refuse to detach time from the motion of things, and leaves it reduced to a simple chronology, *a history*, of these motions.

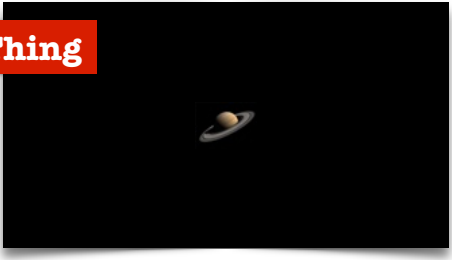

The Thing

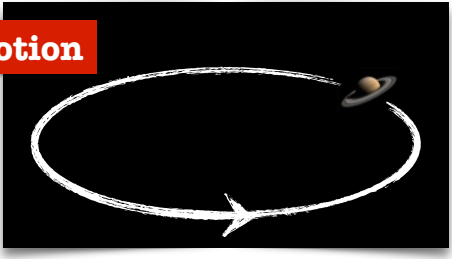

The Motion

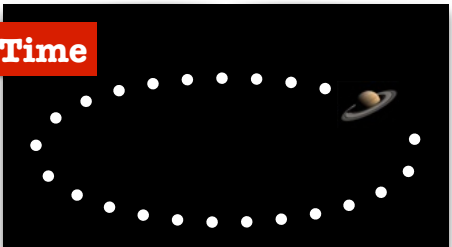

The Time

« **Time is an accident inherent in motion, which is an accident in the thing moved.** »

Maimonides, *The Guide for the Perplexed.*

It is thus that Maimonides is led to refute, according to the principle of minimality, the absolute, abstract (in the etymological sense), character of Aristotelian time, ultimately granting it only the status of an accident of motion. If one wanted to go further, motion itself is only accidental to the thing, it is not its essence. This is why, in particular, speaking of the few moments before the Big Bang, as one sometimes reads, makes no sense and shows a profound ignorance of the fundamental concepts of physics. Unless one imagines that the Big Bang is just an instant in the destruction and reconstruction of a pre-existing thing. Which would quickly lead back to the question of the eternity of the world. Maimonides anticipates this question in chapter VIII : « Indeed, **as soon as you affirm a time before the world,** **you are obliged to admit *eternity***; for time being an accident, to which a *substratum* is necessary, it would follow therefrom that something existed before the existence of this world which now exists, and it is precisely that which we wish to avoid ». But we will not elaborate further on this remark, although it is important and quite modern, insofar as we are only interested in mechanics, *i.e.*, in the description and intelligibility of mechanisms, whether natural (the system of the heavens) or human (pulleys, levers etc.).

However, note that this refutation of Aristotle's absolute time does not concern simultaneity, or the synchronization of clocks, of which Maimonides obviously says nothing and can say nothing. But this does not detract from the relevance of Maimonides' objection, which is a crack concerning the nature of Aristotelian time. It is quite interesting to note that in these centuries where Aristotle's absolute spirit reigned supreme over even the greatest scholars, a doubt regarding the value of his presuppositions on physics persisted in a corner of Spain [13]. But it would take until the end of the 19th century, and Einstein's theory of relativity, to raise this question decisively and provide an operational answer.

**Brief summary regarding time** For *Plato*: Time is generated and born with the heaven. For *Aristotle*: The world is uncreated and eternal, time is detached from motions. For *Maimonides*: It is not proven that the world is eternal, and time cannot be detached from motion, of which it is an accident."'

## 6. Bruno's Objection

We have heard Maimonides' dispute with Aristotle on the nature of time. And if there is a questioning of the validity of Aristotelian time, it does not seem that the conception of space at that time was subject to debate. The idea of a world at the center of which the earth sits, whose immobility defines space, was shared by all scholars until and even after what was called the Copernican revolution, in the sixteenth century. But the new heliocentric vision of the world was not sufficient to challenge the Aristotelian edifice. Certainly, the new cosmology extended space towards the infinite distances and broke the sublunar limit, but it did not at all imply renouncing the notion of the rest of things and thus the Aristotelian notion of space.

---

13. Although Maimonides was born in Cordoba, he wrote his *Guide* in Fostat, Egypt around 1190.

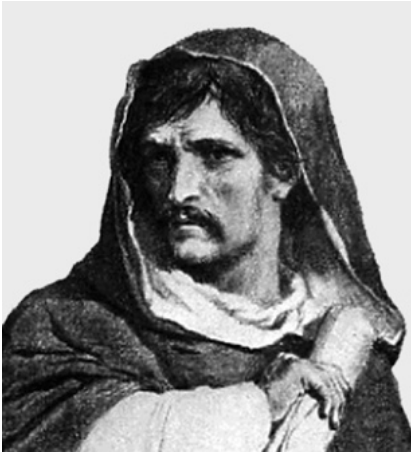

*Giordano Bruno,*
*Nola 1548 — Rome 1600.*

No, the questioning of Aristotelian space is more subtle ; it requires renouncing the very idea of rest, *i.e.*, the idea that we cannot distinguish motion from rest (more precisely, certain motions). It is this idea, more than the place of the earth in the organization of the heavens, that is the hinge of the second major **epistemological rupture** in the history of mechanics, the first being undoubtedly the edition of Aristotle's *Physics*.

It seems that the first (written) elements of this controversy on the nature of Aristotle's Space should be attributed to Giordano Bruno, precisely in his *Ash Wednesday Supper*.[14] It is in this work that he challenges the method by which Aristotle, in the second book of his treatise *On the Heavens*, concludes the immobility of the earth, a method which consists — to put it briefly — in stating : if the earth moves, then upon dropping a stone vertically, it should fall in a different place, further away from its starting point the faster the earth's motion is. To refute this argument, Bruno uses the boat as a metaphor for the Earth. He has Smitho say, in the third dialogue :

> « *SMITHO : The motion that affects the earth must*
> *necessarily change all the relations between straight*

---

14. Giordano Bruno. *The Ash Wednesday Supper*. Éditions de l'éclat, Paris, 2006.

*and oblique lines. There is a difference between the motion of a ship and the motion of what the ship contains. **If it were not so, one could never drag an object in a straight line from one side of the ship to the other** when it is sailing on the open sea; **nor fall back, after a jump, to the exact spot from which one leaped.*** »

The objection is clear; it may not be original, but it had the merit of being written and published. Teofilo clarifies Bruno's thought even further, when he replies :

« *TEOFILO : Thus, everything on the earth moves with the earth* [...] ***A man standing on the bank*** *and throwing a stone straight towards the ship, while it sails down the river,* ***will miss*** *his shot by as much as the speed of the course is greater* [...] *But* ***a man standing on the mast*** *of the ship, however fast the ship is going,* ***will not miss his target : nothing will prevent the stone*** [...] ***from reaching in a straight line*** [...] ***the base of the mast.*** »

In these few sentences is expressed what has since become known as : the ***principle of the relativity of motion***. And Bruno goes a little further in his arguments when he has Teofilo say :

« *TEOFILO : The stone dropped by someone standing on the ship and who is consequently carried along in its motion, is endowed with a virtue that is lacking in the stone dropped by someone standing outside* [...] ***The first stone is endowed with the virtue of the mover who moves with the ship, the other with the virtue of the mover who does not share this motion***. »

The choice of the word "virtue" here is not incidental; it is a tribute physics pays to poetry. Before a concept is flattened into the precise, skeletal language of mathematics, it often begins as a metaphor, a word chosen for its flesh and soul. Bruno's "virtue" is not yet the modern, passive concept of inertia;[15] it is a living quality, an indwelling power that motion imparts to an object. By respecting this original language, we can feel the birth of a new idea before it has been formalized, a profound shift from a purely spatial description of motion towards a more dynamic one.

**Conclusion.** By the experience of the boat, whether he actually experimented with it or not is unknown, Bruno refutes Aristotle's demonstration of the earth's immobility. And if Aristotle's major argument fails for a boat, it indeed has little value concerning the earth : « *Thus, everything on the earth moves with the earth* ».

The impossibility of physically distinguishing between rest and uniform rectilinear motion definitively invalidates the notion of Space as the places of rest. Things may still tend towards their proper places, as Aristotle asserts, but these proper places are increasingly difficult to grasp since one no longer knows if they are not themselves in uniform rectilinear motion. In short, Bruno strikes a serious blow to Aristotle's Space, but it will be up to Galileo to deliver the coup de grâce.

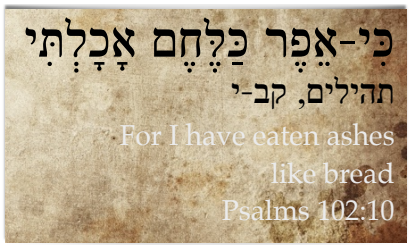

The **Psalm of David 102.10** which inspired the title of *The Ash Wednesday Supper*.

Giordano Bruno was burned alive by the Catholic Church on February 17, 1600, in Campo dei Fiori in Rome, for his insolence and stubbornness in seeking what is true. The biblical reference in the title of his work *The Ash Wednesday Supper* could not better dramatically condense, in so few well-chosen words, the destiny of this noble, tormented soul.

---

15. That is, Galilean inertia, which obviously did not exist at this time. The only concept of inertia that existed, even if not formalized, was the Aristotelian inertia called rest.

# 7.  Galileo's Abandonment of Space

In his *Dialogue Concerning the Two Chief World Systems*[16], Galileo reuses the argument about the experience of the stone dropped from the top of a ship's mast, which we have just discussed through Bruno's text. I do not know, when Galileo speaks of this experience as the *very characteristic experiment*, whether he is referring precisely to Bruno or to a scientific or pre-scientific folklore that might be older. Nevertheless, this experiment is indeed a key to the question of the motion or rest of the earth. And consequently, a key to judging the value of the notion of rest in physics. It is in the *Second Day*

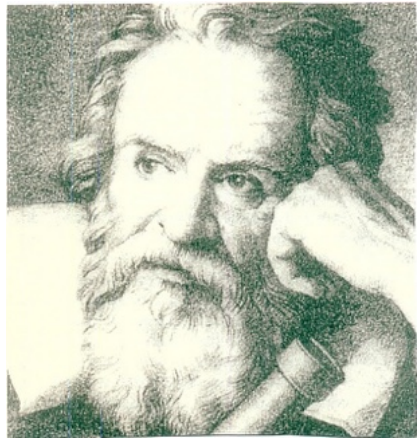

*Galileo Galilei,*
*Pisa 1564, Arcetri 1642.*

(280) that he introduces this debate through his character Simplicio, to whom he attributes Aristotle's argument and has him say :

« *SIMPLICIO : There is also **the very characteristic experience of dropping a stone from the top of a ship's mast** : when the ship is at rest, it falls at the foot of the mast; when the ship is underway, it falls at a distance from the foot equal to the distance the ship has advanced during the time of the stone's fall.* »

Without fully repeating the entire debate surrounding this imaginary experiment, we will side with Galileo's opinion, which he has Salviati express firmly and definitively, still in the *Second Day* (284) of the Dialogue :

> « *SALVIATI : That its authors can present it without*
> *having performed it, you yourself are a good witness :*
> *it is without having performed it that you hold it*
> *for certain, relying on their good faith; it is therefore*

---

16.  Galileo Galilei. *Dialogue Concerning the Two Chief World Systems*, Seuil, collection *Points*, Paris, 1992.

*possible and even necessary that they also, have done the same, I mean that they have relied on their predecessors, without **ever finding someone who has done it**. Let anyone perform it and they will indeed find that experience shows the opposite of what is written : **the stone falls to the same place on the ship, whether it is at rest or moving at any speed**. The same reasoning applying to the ship as to the Earth, if the stone always falls vertically at the foot of the tower, **nothing can be concluded from this regarding the motion or rest** of the Earth.* »

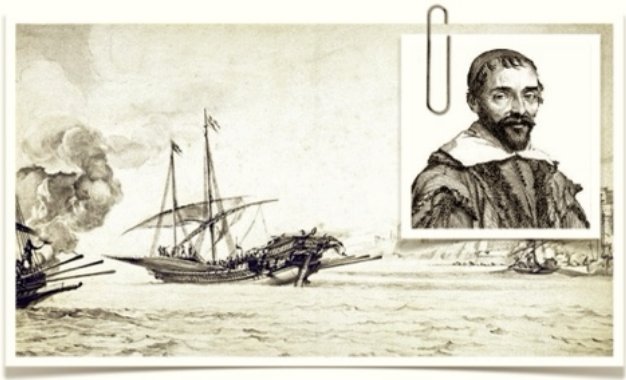

*Gassendi performs the cannonball drop experiment, from the top of the mast of a galley chartered by him, in the port of Marseille in 1641.*

It seems it was Gassendi who first performed this experiment in 1641 on a galley in the port of Marseille, if we are to believe the account of this experiment found in the *Collection of Letters from Sieurs Morin, de La Roche, de Neuré and Gassend*[17]. The result of the experiment indeed confirmed the principle underlying the expectations of Bruno and Galileo : the cannonball dropped from the top of the mast did fall at the base of the mast, neither ahead nor behind, whether the galley was in port or sailing, *with all possible force & speed*, still according to the account. One should read with particular interest Vincent Julien's study, rich in revealing details about this event although — even if one understands the intention — his conclusion is at least debatable,

---

17. Jean-Batiste Morin, Mathurin Neure, Francois de Barancy, Pierre Gassendi. *Recueil de lettres des Sieurs Morin, de la Roche, de Neur, et Gassend etc.* Augustin Courbe, libraire. Paris, 1650.

**AV LECTEVR.**

qui, es rece.

*Monsieur Gassendi ayant esté tousiours tres-curieux de chercher à iustifier par les experiences la verité des speculations que la Philosophie luy propose, & se trouuant à Marseille auec Monseigneur le Comte d'Alais en l'an 1641. fit voir sur vne galere qui sortit exprez en mer par l'ordre de ce Prince, plus illustre par l'amour & la connoissance qu'il a des bonnes choses, que par la grandeur de sa naissance, qu'vne pierre laschée du plus haut du mast, tandis que la galere vogue auec toute la force & la vistesse possible, ne tombe point ailleurs qu'elle ne feroit, si la mesme galere estoit arrestée & immobile; si bien que soit qu'elle aille, ou qu'elle n'aille pas, la pierre tombe tousiours le long du mast à son pié, & de mesme costé. Cette*

*experience faite en presence de Monseigneur le Comte d'Alais, & d'vn grand nombre de personnes qui y assisterent, sembla tenir quelque chose du paradoxe à beaucoup qui ne l'auoient point veuë; ce qui fut cause que Mons' Gassendi composa vn Traité, De motu impresso à motore translato, que nous vismes de luy la mesme année en forme de Lettre escrite à M' du Puy. M' Morin qui auoit fait imprimer quelque temps auparauant son ouurage intitulé, Famosi Problematis de Terræ motu hactenus optata, nunc tandem demonstrata solutio; creut que M' Gassendi n'auoit en autre dessein que d'escrire contre son Liure, pource que dans cette Lettre à M' du Puy, il destruisoit vne des plus fortes raisons que l'on a tousiours opposées au mouuement de la Terre, & que M' Morin employoit pour fondement d'vne de ses principales demonstrations. Ce desplaisir joint à l'ambition qu'il a*

**Account of Gassendi's experiment in the port of Marseille**.
In *Recueil de lettres des Sieurs Morin, de la Roche, de Neur, et Gassend etc.*
By Jean-Batiste Morin, Mathurin Neure, Francois de Barancy, Pierre
Gassendi. Augustin Courbe, libraire. Paris, 1650.

I quote : [18] « *The Marseille galley therefore strictly speaking demonstrates nothing.* » Even if, as some historians note, the experiment does not strictly demonstrate the motion of the Earth, it demonstrates something far more fundamental for our purpose : it confirms Galileo's explicit conclusion that from such an experiment, *nothing can be concluded regarding the motion or rest* of the system [19]. It is a physical vindication of the principle of relativity.

And if there were any remaining doubt on this question of the physical indiscernibility of uniform rectilinear motion, Galileo has Salviati deliver a definitive and luminous speech. It is a comprehensive thought experiment, a complete physics laboratory sealed within the cabin of a 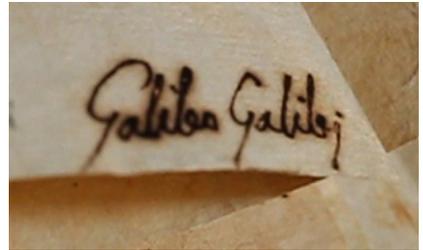 ship, designed to prove that no internal mechanical phenomenon—from the flight of a butterfly to the dripping of water—can ever reveal whether the system is in motion or at rest :

> « *SALVIATI : Shut yourself up with a friend in the largest cabin below deck of a large ship, and take with you flies, butterflies, and other small flying creatures ; also provide yourself with a large container filled with water with small fish ; also hang a small bucket from which water drips into another container with a small opening placed underneath. When the ship is motionless, observe carefully how the small flying creatures move at the same speed in all directions*

*within the cabin, you see the fish swimming indif-
ferently on all sides, the falling drops all enter the
container placed below; if you throw something to
your friend, you don't need to throw harder in one
direction than another when the distances are equal;
if you jump with joined feet, as they say, you will tra-
verse equal spaces in all directions. When you have
carefully observed this, although there is no doubt
that things must happen this way when the ship is
motionless,* **make the ship go at whatever speed you
want; provided the motion is uniform,** *without ro-
cking one way or the other,* **you will not notice the
slightest change in all the effects just mentioned;
none will allow you to tell whether the ship is mo-
ving or immobile** *: when jumping you will traverse
the same distances on the floor as before, and it's
not because the ship is going very fast that you will
make bigger jumps towards the stern than towards
the bow; yet, while you are in the air, the floor be-
neath you moves in the direction opposite to your
jump; if you throw something to your friend, you
won't need more force for him to receive it, whether he
is on the bow or stern side, and you are opposite;* **the
droplets will fall as before** *into the container below
without falling towards the stern, and yet, while the
droplet is in the air, the ship advances several cubits;*
**the fish in their water will not tire more swimming
forward than backward** *in their container, they will
go with the same ease towards the food you have pla-
ced wherever you wish along the edge of the contai-
ner; finally,* **the butterflies and flies will continue
to fly indifferently in all directions,** *you will never
see them take refuge towards the walls on the stern*

> *side as if they were tired of following the rapid course of the ship from which they would have been long separated, since they remain in the air; burn a grain of incense, a little smoke will form which you will see rise upwards and remain there, like a small cloud, without moving one way rather than another. If all these effects correspond, it is because the motion of the ship is common to everything it contains.* »

It is thus that the result of Gassendi's experiment and Salviati's plea put an end to Aristotle's Space, as an imperative category of physics. Indeed, while space as a mathematical construction is legitimate and perfectly conceivable, the fact that it is impossible by a physical experiment to distinguish between uniform rectilinear motion and rest, disqualifies Aristotle's Space as a fundamental category of physics. If there is no longer a privileged state of rest, the very foundation upon which Aristotle's geometry was built—the identification of Space with the set of resting motions—has crumbled. The Group of Aristotle, which was defined by its preservation of this structure, is no longer a valid description of our world. We must now search for a new Inertia Group, one that respects this new, more subtle principle of relativity.

## 8. The Galilean Group

In the passage from the Aristotelian system to the Galilean system, the primitive mathematical structure of spacetime is preserved : it is still a real affine space of dimension 4. We will denote it ET, without this presupposing a particular decomposition into space and time. Points in ET are still called events. We have lost Aristotle's rests [20], but Aristotle's clock $\tau : q \mapsto t$ remains, which associates an event $q \in \mathrm{ET}$ with its date $t \in \mathrm{T}$. Indeed, Galileo's principles do not challenge the absolute universal nature of Time. We will therefore admit :

---

20. Since according to Galileo, *one cannot distinguish rest from uniform rectilinear motion.*

**Statement** *Aristotle's clock* $\tau : \mathrm{ET} \to \mathrm{T}$ *is an affine map.*

*i.e.*, the inverse-images, which we denoted **t** in the Aristotelian case and which we denote $\mathrm{E}_t$ in the Galilean case, of points $t \in \mathrm{T}$ by the clock $\tau$ are affine subspaces of ET. The sheets $\mathrm{E}_t$ are the events that occur simultaneously at instant $t$ ; we can also call these sheets of simultaneous events : the *spaces at instant $t$*, without confounding them with Aristotle's Space, obviously.

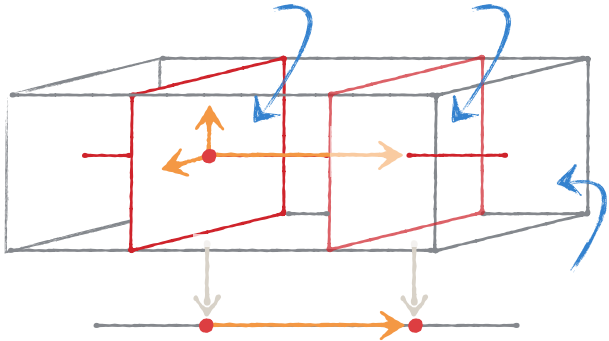

Galilean spacetime

It remains to establish what becomes of distances and durations in this new context. Regarding durations, they are still measured on Time T, since nothing has changed regarding it. One can also obviously measure durations between two events :

**Definition** *The duration* $\delta(q, q')$ *between two events $q$ and $q'$ is defined by*

$$\delta(q, q') = d(\tau(q'), \tau(q)),$$

*where $d$ is the measure of durations introduced previously on time* T.

Regarding distances, it is *a priori* a little more delicate, but the Galilean point of view does not absolutely forbid the measurement of distances between simultaneous events. The fact that there is no absolute reference space in the sense we attribute to Aristotle, does not mean that we are unable to measure the distance between objects in our immediate vicinity. This leads us to consider the existence of a Euclidean metric $d_t$ on each sheet $E_t$. The physical principle of temporal homogeneity demands that these metrics be consistent with each other : a measurement made with a physical ruler should not depend on the absolute time at which it is performed. This implies that the family of distances $d_t$ must be invariant under temporal translations. The following proposition formalizes this intuition.

**Definition** *A temporal translation* in Galilean spacetime is a translation by any vector of the underlying space, transverse to the sheets of simultaneous events.

We then have the following precise statement :

**Proposition.** *There exists a family of Euclidean distances $d_t$ defined on the sheets $E_t$, invariant under temporal translation. i.e., such that for any vector K in ET, transverse to Aristotle's clock, for any $t \in T$ and for any pair of events $q, q' \in E_t$,*

$$d_t(q, q') = d_{t+k}(q + K, q' + K), \qquad (\heartsuit)$$

*where $k$ is the projection of K onto T. This family of distances depends only on the choice of a distance $d$ on an arbitrary origin sheet E, and it is then given by the inverse formula*

$$d_t(q, q') = d(q - K, q' - K), \qquad (\diamondsuit)$$

*for any pair $q$ and $q'$ of events at instant $t$, and for any vector K transverse to the sheets, that translates $q$ and $q'$ into E.*

*Proof.* Let us choose an origin of time, $o \in T$, and denote E the sheet above $o \in T$. Let us then choose a point $o \in E$ that is both the origin of E and ET. Let $d$ be a Euclidean distance on E. Define for any $t \in T$ the distance $d_t$ defined by ($\diamondsuit$) for a certain vector K that translates E to $E_t$. Note that since the translation by K is an isomorphism between E and $E_t$, $d_t$ indeed defines a distance. Let's now show that $d_t$ does not depend on K, as this is what's at stake. Let K$'$ be another vector that translates E to $E_t$, then K$' = $ K $+ \eta$ where $\eta = $ K$' - $K is in the kernel of $\tau$. Then, $d(q - $K$', q' - $K$') = d(q - $K$ - \eta, q' - $K$ - \eta)$, but since $d$ is a Euclidean metric, it is invariant under translation, and $d(q-\text{K}-\eta, q'-\text{K}-\eta) = d(q-\text{K}, q'-\text{K})$. So we have $d(q-\text{K}', q'-\text{K}') = d(q-\text{K}, q'-\text{K}) = d_t(q, q')$. $\qquad\qquad\square$

**In summary** Galilean spacetime ET is equipped with a foliation $\{\text{E}_t\}_{t \in \text{T}}$ into affine sheets of simultaneous events, inverse-images of instants by Aristotle's clock, and a family of Euclidean distances $\{d_t\}_{t \in \text{T}}$ on the sheets, invariant under temporal translations. But that's not all. With the Aristotelian axiom of absolute rest invalidated, the very definition of an 'ideal' or 'unforced' motion must be expanded. This is the axiomatic heart of the Galilean rupture, the move that will generate the new geometry :

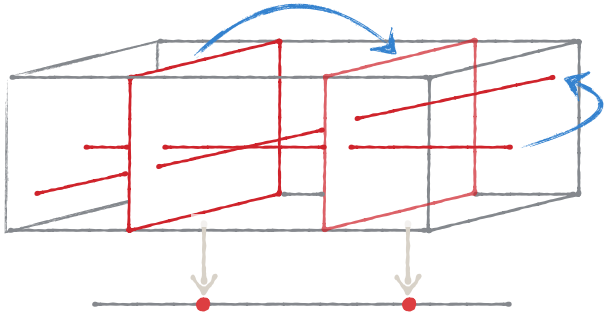

The Galilean Group

**Definition** *Galilean inertial motions* *are uniform rectilinear motions with finite velocity, i.e., the affine lines drawn in* ET *, and transverse to Aristotle's clock.*

All elements are now gathered to introduce the *Galilean Group*, the *Inertia Group* of Galilean Mechanics. It is defined, in parallel to the Aristotelian case, as the set of transformations *g* of ET that satisfy the following conditions :

- ✔ **Preserves Galilean inertial motions** : exchanges uniform rectilinear motions into uniform rectilinear motions.
- ✔ **Preserves Time** : exchanges sheets of simultaneous events.

✔ **Preserves the Euclidean structures** of Time and the sheets of simultaneous events.

**Proposition** *Acting on Galilean spacetime* ET*, the Galilean group is represented by the matrices*

$$\begin{pmatrix} A & B & C \\ 0 & 1 & e \\ 0 & 0 & 1 \end{pmatrix} : \begin{pmatrix} r \\ t \\ 1 \end{pmatrix} \mapsto \begin{pmatrix} Ar + Bt + C \\ \pm t + e \\ 1 \end{pmatrix} \quad where \quad \begin{cases} A \in O(3) \\ B, C \in \mathbf{R}^3 \\ e \in \mathbf{R} \end{cases}.$$

*Proof.* First, we identify Euclidean spacetime with the ordinary product $\mathbf{R}^3 \times \mathbf{R}$, the points of this spacetime are denoted $(r, t)$. It is sufficient to choose a vector L underlying E and transverse to the sheets of simultaneous events for temporal orientation, then choose 3 direction vectors I, J, K of the origin sheet E to construct an affine isomorphism $(r, t) \mapsto o + xI + yJ + zK + tL$, where $o$ is the origin we chose earlier, and $r = (x, y, z)$. By the Gram-Schmidt orthonormalization process, it is possible to choose base vectors such that the Euclidean distance on $E \sim \mathbf{R}^3$ and the duration measure on $T \sim \mathbf{R}$ are given by formulas (♣) and (♠). Let $g$ be a transformation of $\mathbf{R}^3 \times \mathbf{R}$ satisfying the Galilean conditions.

Next, since $g$ maps every affine line to another affine line, it is therefore an affine map. Denoting $X = (r, t)$, we can write the transformation matricially :

$$g(X) = \begin{pmatrix} A & B \\ b & a \end{pmatrix} \begin{pmatrix} r \\ t \end{pmatrix} + \begin{pmatrix} C \\ e \end{pmatrix} = \begin{pmatrix} Ar + Bt + C \\ br + at + e \end{pmatrix}$$

But since $g$ exchanges the sheets of simultaneous events and the transformation on T preserves durations, for the same reasons as in the Aristotelian case, we have $b = 0$ and $a = \pm 1$. Which finally gives :

$$g(X) = \begin{pmatrix} Ar + Bt + C \\ \pm t + e \end{pmatrix}, \quad \text{with} \quad B, C \in \mathbf{R}^3, \ e \in \mathbf{R}.$$

It remains to translate the conservation of the Euclidean structure on the sheets of simultaneous events, which gives the last element $A \in O(3)$. □

**Note however** that it is customary to define the Galilean group rather as the connected component containing the identity, *i.e.*, with $A \in SO(3)$ and $t \mapsto t + e$.

**Galilean boosts :** As we can see, the Galilean group is a supergroup of the Aristotelian group. The matrix block B that is added is responsible for what are called *Galilean boosts*, or the act of *imparting motion*. These transformations map Aristotelian inertial motions [21] to Galilean inertial motions : the constant motion $t \mapsto r$ is mapped to the boosted motion $t \mapsto r + B\,t$. This is how this group captures Galileo's principle of the impossibility of rest.

**Spacetime is still a homogeneous space** of the inertia group of Galilean mechanics. It is still the primary domain of Galilean geometry.

**What then remains of Aristotle's Space?** Aristotle's rests no longer constituting an invariant space of the inertia group, we witness the disappearance of Space, at least as a singular class of motions. However, the mathematical impossibility of an acceptable substitute for this space, for which we all, despite everything, have a natural sensitivity, remains to be expressed. Since Galilean spacetime has preserved its nature as an affine space, it is its decomposition $ET \simeq E \times T$ that becomes problematic. The Aristotelian clock, the natural projection $\tau\colon ET \to T$, remains legitimate, however, since, by construction, the Galilean group exchanges the sheets of simultaneous events. Thus, it is the projection

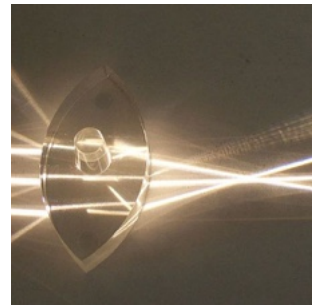

**Light** Lines that are not transverse to Aristotle's clock are entirely contained within the sheets of simultaneous events. These are instantaneous motions and they represent *light rays* in Galilean geometrical optics.

from ET onto a hypothetical space E of dimension 3 that ultimately summarizes the question of space in Galilean mechanics. Does there then exist a projection, from spacetime ET onto a Euclidean space E of dimension 3, that could play the role of Space in this new context? To be consistent, this projection should then intertwine the actions of the Galilean group and the Euclidean group.

It is the Galilean boosts that prevent the existence of such a projection, as shown by the following proposition.

---

21. This assumes an arbitrary identification of temporal sheets at each instant with one of them at a given instant, but this is not crucial.

**THEOREM** *There exists no submersion*[22] *from* ET *onto* E *that intertwines the actions of the Galilean group with the Euclidean group.*

*Proof.* Let $\pi\colon \mathrm{ET} \to \mathrm{E}$ be a submersion that intertwines the actions of the Galilean group and the Euclidean group. *i.e.*, $\pi(gx) = h(g)\pi(x)$ for any element $g$ of the Galilean group and any $x \in \mathrm{ET}$. By equivariance, the fibers of $\pi$ are exchanged by the Galilean group, *i.e.* : if $\pi(x) = \pi(x')$ then $\pi(gx) = \pi(gx')$. In particular, if $V_x$ denotes the kernel of the tangent linear map $\mathrm{D}(\pi)(x)$, then $V_{gx} = gV_x$. Indeed, from $\pi(gx) = h(g)\pi(x)$ we obtain infinitesimally $\mathrm{D}(\pi)(gx)(g\xi) = h(g)\mathrm{D}(\pi)(x)(\xi)$, for any point $x \in \mathrm{ET}$ and any vector $\xi$ in $\mathbf{R}^4 \times \{0\}$. Thus, if $\xi \in \ker(\mathrm{D}(\pi)(x))$, then $g\xi \in \ker(\mathrm{D}(\pi)(gx))$, and since $g$ is invertible, we have equality. Let's place ourselves at the origin $\mathbf{0}$ of ET, represented matricially by the vector $(0_3, 0, 1)$ in our affine model. Consider then the stabilizer S of $\mathbf{0}$, *i.e.* the subgroup of elements $g$ of the Galilean group such that $g\mathbf{0} = \mathbf{0}$. According to the preceding, we should have $gV_{\mathbf{0}} = V_{\mathbf{0}}$ for all $g \in \mathrm{S}$. Let $g \in \mathrm{S}$, then :

$$g = \begin{pmatrix} \mathrm{A} & \mathrm{B} & 0 \\ 0 & 1 & 0 \\ 0 & 0 & 1 \end{pmatrix}, \text{ with } \mathrm{A} \in \mathrm{SO}(3) \text{ and } \mathrm{B} \in \mathbf{R}^3.$$

Its infinitesimal action on a vector $\xi = (v, \varepsilon, 0) \in \mathbf{R}^4 \times \{0\}$ tangent to ET at the origin is the linear action

$$g\xi = \begin{pmatrix} \mathrm{A}v + \varepsilon\mathrm{B} \\ \varepsilon \\ 0 \end{pmatrix}, \quad \text{i.e.} \quad \begin{pmatrix} \mathrm{A} & \mathrm{B} \\ 0 & 1 \end{pmatrix}\begin{pmatrix} v \\ \varepsilon \end{pmatrix} = \begin{pmatrix} \mathrm{A}v + \varepsilon\mathrm{B} \\ \varepsilon \end{pmatrix}.$$

We have assumed that the projection $\pi$ is a submersion onto E, *i.e.*, that the kernel of its tangent linear map is one-dimensional everywhere. And thus that the action of S has a one-dimensional invariant subspace; however, this is not the case, essentially because SO(3) does not have one. Therefore, the desired intertwining does not exist. □

In other words, **there is no Galilean replacement for Aristotle's Space**. From the Galilean perspective, there is no reasonable way to give meaning to the classic and misleading sentence from last century's science fiction novels : « And the rocket stopped in space »… Rockets no longer stop in space, since 1632 [23].

---

22. Technically, a surjective map whose tangent linear map is itself maximal at every point, *i.e.*, surjective and non-singular.

**Is Galilean space then six-dimensional?** With Galileo, we lost Aristotle's three-dimensional Euclidean Space, without even the possibility of finding an equivalent. Out of respect for Galileo's new principles on motion, we have therefore renounced Aristotle's inertial motions as structuring Space. It is then legitimate, but perhaps not entirely intuitive, to define Galilean space as the set of all Galilean inertial motions, *i.e.*, all uniform rectilinear motions. We immediately notice that this is a 6-dimensional space, since a uniform rectilinear motion is uniquely defined by an initial position $r$ and a velocity $v$.

We will see in the second part of this book that this six-dimensional space of inertial motions is the true "space" of Galilean mechanics. It is no longer endowed with a Euclidean structure, as Aristotle would have wished, but with a canonical **symplectic structure**.

This structure, first uncovered by Lagrange [24], is the geometric key to dynamics and will be the central object of our investigation. It is also at the origin of an entire field of mathematics very rightly called « Symplectic Geometry ».

## 9. Return to Time, Einstein…

We have seen how Maimonides had already questioned Aristotle's Time, to which the motion of everything is supposed to relate. His argument, even though remarkably insightful and premonitory,

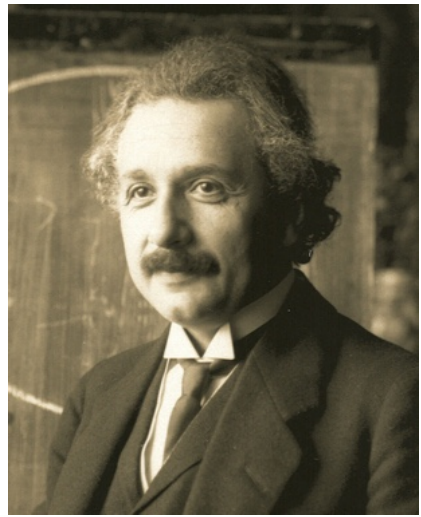

*Albert Einstein*
*Ulm 1879, Princeton 1955*

---

23. As J.-M. Souriau systematically and emphatically noted.

24. Patrick Iglesias. *The origins of symplectic calculus in Lagrange*. L'enseignement mathématique, volume 44, pages 257—277, 1998.

was mostly an argument by default : There is no time outside of motion. If there is an eternal time that can serve as a universal reference, then there is a body in act whose motion is eternal. This idea of eternity, Maimonides could not accept without debate, and finding some weakness in Aristotle's argumentation, he preferred to adhere, with Plato, to the idea that there was beginning and that time was created with the world (the Heavens as said Plato). Remarkably, what is known today through experience : the 3°K cosmic microwave background — the Big Bang — the original explosion, and the 15 billion years of existence of our universe, all this gives *a posteriori* reason to Maimonides, whose respect for the principle of minimality of axioms is thus rewarded [25].

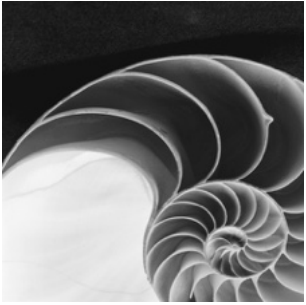

**The Big Bang** is a relativistic spacetime singularity in cosmology — a 4-dimensional Lorentzian manifold. A singularity can be defined as a boundary at infinity (obtained by exhausting compact sets) towards which geodesics converge in finite time. The boundary must not be regularizable, *i.e.*, it must not be possible to eliminate it by embedding it in another spacetime. *The image intended to suggest the idea of a singularity is an X-ray photograph of a snail shell.*

The refutation of Aristotle's time was ultimately organized on two fronts : on the one hand, we cannot trace back through the Big Bang, and even the most reluctant have accepted the idea of a beginning 15 billion years ago. Therefore, there is no body in act whose eternal motion could serve as a reference. But even more serious than that, all our previous constructions implicitly assumed that we could legitimately [26] foliate spacetime into sheets of simultaneous events. We will see that even this was asking too much of nature.

It took until the 19th century, with the questions raised by the variance of Maxwell's equations of electromagnetism, and the contradictions raised by certain experiments, the most famous of which is the Michelson-Morley experiment, for the beautiful Aristotelian temporal edifice to collapse, much to the chagrin of physicists of the time,

---

25. One could invoke here the famous *Ockham's Razor* if it were not slightly anachronistic ; Maimonides died in 1204 and William of Ockham was born in 1287.

26. Legitimate, from the point of view of the laws of nature.

and for them, after being forced to renounce Space, to be led to renounce Time.

**The Michelson-Morley experiment.** Despite Galileo, the idea of an absolute space never ceased to haunt the minds of physicists. This unlikely quest was renewed with the wave theory of light and Maxwell's equations, at the end of the nineteenth century. Since light was ultimately described as a wave, there must exist a subtle medium whose elastic properties it manifested. This medium must permeate all physical bodies, fill space, and light must propagate through it in all directions with the same speed $c$. A new absolute space thus came to replace the defunct Aristotelian Space : the ***luminiferous Aether***. In a sense, this was the last great attempt to find a physical substance for Aristotle's principles—a universal, immobile medium whose properties would define an absolute frame of reference for all motion and, through it, a universal Time.

Michelson decided to measure the speed of the Earth relative to this luminiferous aether using an experiment that now bears his name, and which consists of splitting a monochromatic light ray and making it travel two orthogonal paths before recombining them, and measuring the shift between the two beams. The experimental apparatus is diagrammed in the sketch below. The shift at the end of the path between the longitudinal ray (red in the figure) and the transverse ray (blue in the figure) is now classic, even though Michelson had initially given an incorrect version. He had indeed expressed the distance traveled

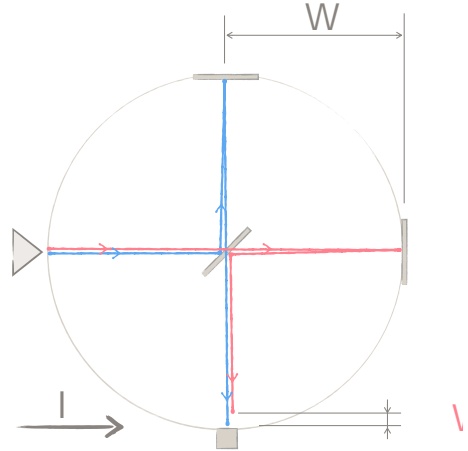

**The Michelson Interferometer** A monochromatic beam is sent onto a semi-transparent central mirror where it is split into two orthogonal beams which are then reflected back by the mirrors towards the central mirror to be recombined and directed onto an interference detector which is supposed to register the phase difference between the two beams.

by the transverse ray as if it were attached to the Earth, contrary to the very hypotheses that the experiment was supposed to test. While this error (corrected by Alfred Potier in 1882) ultimately had no impact on the conclusion, since in any case the results of the experiment were negative, it nevertheless shows how fragile the hypothesis of this immobile aether was, which carried the propagation of light waves while everything else moved independently around or through it.

$$d_t = \frac{2d}{\sqrt{1-\beta^2}} \quad \text{and} \quad d_l = \frac{2d}{1-\beta^2} \quad \text{with} \quad \beta = \frac{v}{c}.$$

The difference $\Delta d$ between the two optical paths $d_t$ and $d_l$ is therefore of the order of $(v/c)^2$, which is infinitesimally small [27] and difficult to measure as such, practically impossible.

Michelson then had the idea of processing this optical path difference with an interferometer placed at the end of the path. Since light is wave-like, it generates interference patterns, and when related to frequencies, this difference becomes sensitive to measuring instruments and significant.

Albert Michelson developed his interferometer in Potsdam where he conducted his first experiment in 1881. The original idea of measuring the Earth's velocity relative to the luminiferous aether, by comparing decoupled optical paths, is probably due to James Clerk Maxwell, the founder of modern electromagnetism. Here is what Robert S. Shankland says about this [28], in his address to the *Michelson Colloquium* held in Potsdam in 1981.

> « *However, there is a very direct connection between Michelson's work and a letter received a year earlier on March 19, 1879 from James Clerk Maxwell addressed to David Peck Todd at the National Almanac Office where Michelson worked with Newcomb.*

---

27. This speed was thought to be small, of the order of the Earth's rotational speed relative to the sun, around a few tens of kilometers per second.

28. Robert S. Shankland. *Michelson in Potsdam*. Astronomische Nachrichten, vol. 303, pp.3–5, Weinheim, 1982.

> *Maxwell had enquired whether observations on the satellites of Jupiter could be made with sufficient precision to extend the velocity of light determinations of Roemer by making them at times near dawn and dusk to reveal changes in the speed of light caused by the earth's motion through the aether. Maxwell explained that his proposed method would only require accuracy to the order ($v/c$), the ratio of the earth's speed to that of light, whereas all purely terrestrial aether-drift experiments demanded the much greater sensitivity of ($v/c$)$^2$, and that no apparatus then existed capable of this severe requirement. **But this assertion of Maxwell's was the challenge that the young Michelson accepted and which led to the invention of his interferometer!*** »

Despite Michelson's remarkable efforts [29], his experiment never revealed interference, although it has been repeated several times since 1881, by himself or by other physicists. This failure sanctioned the attempt to replace Aristotle's Space with a tangible, immobile, and eternal luminiferous Aether [30].

With the « failure » of this experiment, the question of the speed of light, of its measurement in different embedded systems [31]remains open. It is by trying to resolve this contradiction that Hendrik Lorentz,

---

29. Efforts rewarded with the Nobel Prize awarded to him in 1907.

30. A positive result from this experiment would challenge Galileo's conclusions, by providing an electro-mechanical means of determining whether a body is at rest or in uniform rectilinear motion. Meanwhile, the scientific community has accepted, sometimes reluctantly, given how deeply the idea of a Space is ingrained in minds, the abandonment of the aether as a subtle but tangible and resting luminiferous substrate.

31. The notion of *embedded system* is the active version of the notion of *frame of reference*.

Henri Poincaré [32], and Albert Einstein established the principles of the theory of relativity. This theory, which involves a disconcerting mixture of Space and Time, never ceases to surprise, even after more than a century of history.

Let's return to the Michelson experiment. Hendrik Lorentz was the first to propose a solution [33], even if it may seem strange to us *a priori* : he suggests, to compensate for the delay of the longitudinal beam compared to the transverse beam, to state that, under the effect of a mysterious force exerted by the aether on matter, the arm of the interferometer, which moves longitudinally in the direction of the Earth's motion, contracts by the appropriate ratio. Thus, denoting $d'_l$ the length of the contracted arm, to have equality of travel times $d_t/c = d'_l/c$, we need, taking into account the expressions and notations introduced above :

$$\frac{2d'}{1-\beta^2} = \frac{2d}{\sqrt{1-\beta^2}} \quad \text{and thus :} \quad d' = d\sqrt{1-\beta^2}.$$

This is what has since been called Lorentz's *length contraction law*. We see in particular that the speed of light $c$ becomes in this way an insurmountable limit [34]. Regarding the legitimacy of such an hypothesis, on the contraction of moving bodies, here is exactly what Lorentz says about it [35]

---

32. Henri Poincaré, *On the Dynamics of the Electron* (Sur la dynamique de l'électron). C.R Acad. Sci. Paris. Presented June 5, 1905. vol. 140, pp. 1504–1508, Paris (1905). Published in : Rend. Circ. Matem. Palermo 21, pp. 17–176. Palermo (1906) (Received July 23, 1905).

33. Hendrik Lorentz. *The relative motion of the Earth and the Aether* (De relatieve beweging van de aarde en den aether). Zittingsverlag Akad. v. Wet., vol. 1, p. 74, Amsterdam (1892).

34. Unless one thinks that by exceeding the speed of light, bodies which have shrunk to the point of having no thickness, develop in an imaginary dimension...

35. In Lorentz's language : « *Waardoor toch worden de grootte en de gedaante van een vast lichaam bepaald ? Klaarblijkelijk door de intensiteit der moleculaire krachten ; elke oorzaak die deze wijzigde zou ook op den vorm en de afmetingen invloed hebben. Nu mogen wij tegenwoordig wel aannemen dat electrische en magnetische*

> « *Indeed, what determines the size and shape of a solid body? Apparently the intensity of the molecular forces; any cause that could modify them, could also modify the form and dimensions. We may now assume that electric and magnetic forces act through the intervention of the aether. It is not unnatural to assume the same for molecular forces, but then it can make a difference, depending on whether the line connecting two material particles, which move together through the aether, moves parallel to the direction of motion or is perpendicular to it.* »

**Length contraction.** For Lorentz, length contraction is therefore not an optical illusion, or a theoretical convenience to circumvent an embarrassing result, but the consequence of a real physical phenomenon : the *friction*, the *pushing*, of the luminiferous aether on moving matter. And even if one can be satisfied with this explanation, the principle of length contraction alone does not resolve all the questions raised by the hypothesis of the luminiferous aether. In particular, the Trouton-Noble experiment [36] poses a new challenge in this regard. Trouton had undertaken to detect the rotational torque, acting on a flat capacitor, induced by the Earth's motion relative to the aether. The electric charges on the plates of a flat capacitor, placed on the laboratory bench, are in motion relative to the aether, and are therefore the source of a magnetic field. If this is the case, the capacitor is subjected to an electromagnetic force, partially due to the aether wind, and which, according to the laws of electromagnetism, should align the plates of the capacitor perpendicularly to the direction of the Earth's motion relative to the aether. The experiment was carefully conducted between 1901 and 1903, without yielding the expected results. No aether wind blew strongly enough to orient the capacitor.

---

*krachten door tusschenkomst van den aether werken. Het is niet onnatuurlijk hetzelfde voor de molekulaire krachten te onderstellen, maar dan kan 't een verschil maken of de verbindingslijn van twee stofdeeltjes, die zich te zamen door den aether verschuiven, evenwijdig aan de bewegingsrichting loopt of loodrecht daarop staat.* »

36. Frederik T. Trouton and Henry R. Noble. *The mechanical forces acting on a charged electric condenser moving through space.* Phil. Trans. Royal Soc. A 202, pp. 165—181, London (1903).

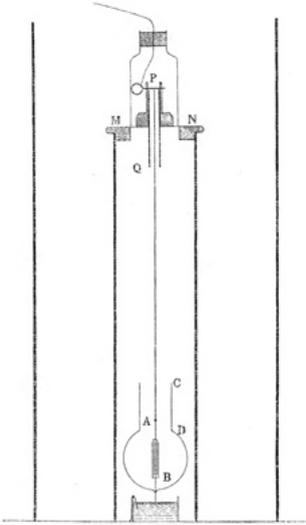

**The Trouton-Noble Capacitor.**
« *PA, the suspension, is a phosphor bronze strip 37 centims. long, the finest that could be obtained. This was soldered at its lower end A to a copper cap, fixed to the condenser protectmg the projecting tin-foil tags, and making contact with them by means of fusible metal. The upper end of the suspension, P, was wound on a small windlass, which was insulated by a mica plate fixed to an annular wooden ring MN, forming the lid to the inner zinc vessel … »*

**Lorentz transformations.** One quickly realizes that the principle of length contraction, which could explain the absence of interference fringes in the Michelson-Morley experiment, is of no help in explaining the negative result of the Trouton-Noble experiment. It is to this new challenge that Lorentz dedicates his article [37] *Electromagnetic phenomena in a system moving with any velocity smaller than that of light*, published in 1904. Poincaré clearly summarized the author's intentions in his article *On the Dynamics of the Electron* of 1905 (*op. cit.*) which he devoted to it, and it is futile to attempt to paraphrase him :

« *Lorentz's idea can be summarized as follows : if one can, without any of the apparent phenomena being modified, impart a common translation to the entire system, it is because the equations of an electromagnetic medium are not altered by certain transformations, which we will call Lorentz transformations; two systems, one at rest, the other in translation, thus become exact images of each other.* »

37. Hendrik A. Lorentz. *Electromagnetic phenomena in a system moving with any velocity smaller than that of light.* Proceedings of the Academy of Amsterdam, vol. 6, 1903-1904, pp. 809-831. Amsterdam (1904).

This is what Lorentz proposes to establish : the formulas for coordinate transformation from the resting ae­ther system to a uniform rectilinear moving system such that the formal structure[38] of Maxwell's equations is preserved un­der the transformation. In other words — if that makes sense without further precision — the equations for the electroma­gnetic field and forces must ultimately have the same « form » in the moving embedded system as in the immobile aether.

These transformations, which Poincaré[39] would call *Lorentz transformations,*[40] must then allow describing how the ob­jects of the electromagnetic system : va­lues of the electric and magnetic fields, charge density, potentials, distances, du­rations etc. transform.

This series of transformations heralds a new mechanics whose laws mix, against all intuition, space and time measurements, and thereby put an end to the idea of a universal reference Time. As we will see, there will no longer be a way in the new mechanics to coherently and permanently synchronize two watches that follow signi­ficantly different motions. This is probably the most unexpected consequence of this

$$div\ \mathfrak{d} = \varrho \quad , \quad div\ \mathfrak{h} = 0,$$

$$rot\ \mathfrak{h} = \frac{1}{c}\,(\dot{\mathfrak{d}} + \varrho\,\mathfrak{v}),$$

$$rot\ \mathfrak{d} = -\frac{1}{c}\,\dot{\mathfrak{h}},$$

$$\mathfrak{f} = \mathfrak{d} + \frac{1}{c}\,[\mathfrak{v}.\,\mathfrak{h}].$$

**Maxwell's Equations.**

The letter $\mathfrak{d}$ represents the « *dielectric displa­cement in the aether* », and the letter $\mathfrak{h}$ the « *magnetic force* ». The dot above represents the time derivative. The letter $\varrho$ represents the volume charge density of an electron, and $\mathfrak{v}$ the velocity of the electron ; the letter $\mathfrak{f}$ the electric force per unit charge that the aether exerts on a volume element of the electron, and the bracket $[\mathfrak{v}.\mathfrak{h}]$ represents the cross pro­duct of the two vectors.

38. The notion of *formal structure* of a system of equations is quite vague, even more so at this time when the foundations of differential geometry were lacking. But the goal pursued is clearly formulated in Poincaré's article.

39. We will see that the transformations proposed by Lorentz had some weak­nesses which Poincaré fortunately corrected.

40. While Lorentz wrote down the coordinate transformations, it was Poincaré who first recognized their group structure and named them after Lorentz. The full affine group, including translations, is now universally known as the Poincaré group.

adventure. If the disappearance of space
had been announced long ago, and was ultimately only confirmed by
the unsuccessful attempts to immerse us in an ethereal fog of light,
the same is not true for Time, which had resisted until then.

Lorentz transformations are precisely defined by Poincaré by formulas (3) of his article. They are generated on the one hand by translations

$$(x, y, z, t) \mapsto (x + C_x, y + C_y, z + C_z, t + C_t),$$

and on the other hand, by transformations $(x, y, z, t) \mapsto (x', y', z', t')$ of the type

$$\begin{cases} x' = k(x + \varepsilon t) & y' = y \\ t' = k(t + \varepsilon x) & z' = z \end{cases}, \quad \text{with} \quad k = \frac{1}{\sqrt{1 - \varepsilon^2}},$$

*preceded and followed by a suitable rotation,* according to Poincaré's own words [41]. As can be seen, the transformation mixes space and time coordinates in an unusual and surprising way, the consequences of which will be measured later. It is these transformations that allow Maxwell's equations to keep their formal appearance, at rest or in uniform rectilinear motion.

One can see, by comparing the texts of Lorentz and Poincaré below, that the transformations they each propose do not exactly coincide. Poincaré also mentions an error by Lorentz in the transformation of charge density [42] which is certainly the result of this shift in the expression of the transformation itself. Here is precisely what Poincaré writes in his introduction by way of caution : « *the results I have obtained are in agreement with those of Mr. Lorentz on all important points ; I have only been led to modify and complete them in a few details* ».

---

41. On the Dynamics of the Electron (Rendiconti etc.), *op. cit.* p. 146.

42. These densities are given by formulas No. 7 in Lorentz's text $\rho' = \rho/(k l^3)$, and No. 4 in Poincaré's text $\rho' = k\rho(1 + \varepsilon\xi)/l^3$.

e shall further transform these formulae by a change of Putting

$$\frac{c^2}{c^2 - w^2} = k^2, \quad \cdots \cdots \quad (8)$$

and understanding by $l$ another numerical quantity, to be determined further on, I take as new independent variables

$$x' = k\,l\,x \quad, \quad y' = l\,y \quad, \quad z' = l\,z, \quad \cdots \quad (4)$$

$$t' = \frac{l}{k}\,t - k\,l\,\frac{w}{c^2}\,x, \quad \cdots \quad (5)$$

Ces équations sont susceptibles d'une transformation remarquable découverte par LoRENTZ et qui doit son intérêt à ce qu'elle explique pourquoi aucune expérience n'est susceptible de nous faire connaître le mouvement absolu de l'univers. Posons :

(3) $\qquad x' = k\,l(x + \varepsilon t), \qquad t' = k\,l(t + \varepsilon x), \qquad y' = l\,y, \qquad z' = l\,z,$

$l$ et $\varepsilon$ étant deux constantes quelconques, et étant

$$k = \frac{1}{\sqrt{1 - \varepsilon^2}}.$$

As for the transformations themselves and their strangeness, note that for Lorentz, the variable $t$ represents the true time, « *the true time* », this is the time attached to the immobile luminiferous aether. And concerning the variable $t'$ : « *The variable $t'$ may be called the "local time"* ». The schema for Lorentz is therefore still Aristotelian : the luminiferous aether defines an electromagnetic Space, immobile and eternal, which replaces Aristotle's mechanical space attached to the earth. The same aether also defines an electromagnetic Time, which replaces the mechanical time of Aristotle's celestial sphere,

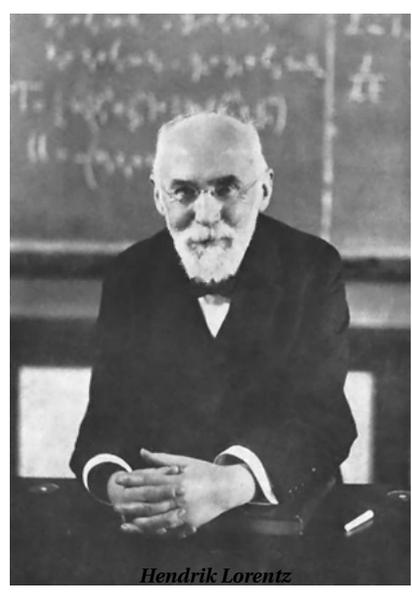

*Hendrik Lorentz*
*Arnhem 1853, Haarlem 1928*

and to which the various local times of embedded systems will be compared.

With the study of these transformations, Poincaré completes the work begun by Lorentz : If no experiment has so far been able to reveal the Earth's motion relative to the aether, it is because uniform rectilinear motions are transparent to the laws governing electromagnetism.

**Einstein's Relativity.** In his historic article « On the Electrodynamics of Moving Bodies » [43] which lays the foundations of his *Theory of Relativity*, Albert Einstein seems more concerned with the lack of symmetry in the treatment of magnets and conductors, regarding electromagnetic forces, than with the contradictions raised by Michelson's or Trouton's experiments, as were Lorentz and Poincaré. Here is what he says in his introduction [44] :

> « *It is known that Maxwell's electrodynamics, as usually understood at present, when applied to moving bodies, **leads to asymmetries which do not appear to be inherent in the phenomena**. Take, for example, the reciprocal action of a magnet and a conductor. The observable phenomenon here depends only on the relative motion of the conductor and the magnet, whereas the customary conception makes a sharp distinction between the two cases where the one or the other of these bodies is the moving one.* »

---

43. Albert Einstein. *Zur elektrodynamik bewegter körper*. Annalen der Physik 17 (10), pp. 891–921. Leipzig 1905 — Submitted June 30 and published September 26, 1905.

44. *Daß die Elektrodynamik Maxwells – wie dieselbe gegenwärtig aufgefaßt zu werden pflegt – in ihrer Anwendung auf bewegte Körper zu Asymmetrien führt, welche den Phänomenen nicht anzuhaften scheinen, ist bekannt. Man denke z. B. an die elektrodynamische Wechselwirkung zwischen einem Magneten und einem Leiter. Das beobachtbare Phänomen hängt hier nur ab von der Relativbewegung von Leiter und Magnet, während nach der üblichen Auffassung die beiden Fälle, daß der eine oder der andere dieser Körper der bewegte sei, streng voneinander zu trennen sind.*

Even if this lack of symmetry obviously bothers and motivates Einstein, he is not unaware of the unfortunate results of the experiments we have discussed so far, and he incorporates them into his exposition to establish his conclusions [45] :

> « *Examples of this type, and the unsuccessful attempts to detect a motion of the Earth relative to the "light medium", suggest that* **electromagnetic phenomena as well as mechanical phenomena possess no properties corresponding to the idea of absolute rest**. *They suggest rather that the same laws of electrodynamics and optics are valid for all coordinate systems for which the equations of mechanics hold, as has already been proven for magnitudes of the first order.* »

These defects, in the expression of electrodynamics and the question about the Earth's velocity relative to the aether, Einstein will correct based on two premises that he will elevate to the rank of principles, and which are : firstly — as for Lorentz and Poincaré — the permanence of the form of Maxwell's equations ; and secondly, the constancy of the speed of light in embedded systems [46]

$$U = \frac{v + w}{1 + \dfrac{v\,w}{V^2}}.$$

**Einstein's Velocity Addition Law**
« *It is remarkable that v and w enter into the expression for the resultant velocity in a symmetric manner.* »

---

45.  *Beispiele ähnlicher Art, sowie die mißlungenen Versuche, eine Bewegung der Erde relativ zum „Lichtmedium" zu konstatieren, führen zu der Vermutung, daß dem Begriffe der absoluten Ruhe nicht nur in der Mechanik, sondern auch in der Elektrodynamik keine Eigenschaften der Erscheinungen entsprechen, sondern daß vielmehr für alle Koordinatensysteme, für welche die mechanischen Gleichungen gelten, auch die gleichen elektrodynamischen und optischen Gesetze gelten, wie dies für die Größen erster Ordnung bereits erwiesen ist.*

46.  *Das Licht im leeren Raume stets mit einer bestimmten, vom Bewegungszustande des emitterenden Körpers unabhängig Geschwindigkeit V fortpflanze.*

> « *Light in empty space always propagates with a definite velocity* V *which is independent of the state of motion of the emitting body.* »

It is this principle that is the basis of Einstein's velocity addition formula, which expresses part of the Lorentz transformations. He will also face the appearance of a **proper time**. He will give it a real physical interpretation, radically different from that of Lorentz or Poincaré, and which will have serious consequences that have since been experimentally verified. We read in his text (*op. cit.* p. 905) [47] :

> « *From this we conclude that a balance-wheel clock situated at the Earth's equator must go slower by a very small amount than a precisely similar clock situated at one of the Earth's poles under otherwise identical conditions.* »

$$\tau = t \sqrt{1 - \left( \frac{v}{V} \right)^2}$$

**Einstein's Proper Time**

« *The clock (when viewed from the stationary system) goes slower each second by* $(1 - \sqrt{1 - (v/V)^2})$ *seconds, or—neglecting magnitudes of the fourth and higher order—by* $\frac{1}{2}(v/V)^2$ *seconds.* »

From the permanence of Maxwell's equations in inertial frames, he deduces an important conceptual (and practical) consequence which he formulates without ambiguity (*op. cit.* p. 892) [48] :

« *The introduction of a "luminiferous ether" will prove to be superfluous, since, according to the view to be developed, neither an "absolutely stationary space" endowed with special properties, nor etc.* »

---

47. *Man schlieBt daraus, daß eine am Erdaquator befindliche Unruhuhr urn einen sehr kleinen Betrag langsamer laufen muß als eine genau gleich beschaffene, sonst gleichen Bedingungen unterworfene, an einem Erdpole befindliche Uhr.*

48. *Die Einfuhrung eines „Lichtathers" wird sich insofern als uberflüssig erweisen, als nach der zu entwickelnden Auffassung weder ein mit besonderen Eigenschaften ausgestntteter „absolut ruhender Raum" eingeführt, noch etc.*

Thus, Einstein's principles of relativity are part of the Galilean logic [49] on the impossibility of distinguishing rest from inertial motion. In this, he distinguishes himself from the *neo-Aristotelian* viewpoint of Lorentz and Poincaré, who continue to base their principles on an immobile aether, an absolute reference for all motion.

NOTE. Without entering into a futile debate about the authorship of the theory of relativity, which is still not free from questionable ideological motivations [50], one can observe, by reading the founding articles, a gradual maturation of the very concept of the « Principle of Relativity ». Its premises can be found in Lorentz's article *Electromagnetic phenomena in a system …*, published in 1903, where he writes [51] (*op.cit.* p. 819) :

> « *In the second place I shall suppose that the forces between uncharged particles, as well as those between such particles and electrons, are influenced by a translation in quite the same way as the electric forces in an electrostatic system.* »

Indeed, how to explain that only charged particles are sensitive to the pressure of the aether wind ; one should then observe an inhomogeneous behavior of matter depending on its constitution, and this is not the case. It is clear that these transformations — which Lorentz has just established — must apply to matter, indiscriminately. But this « hypothesis » is not yet the principle of relativity in all its power, rather a kind of compromise between theory and observation.

---

49. Even if technically Einstein's relativity is an epistemological break of the first order compared to Galilean mechanics, it confirms this specific aspect of the indifference of inertial motion.

50. For an overview of the debate, one may profitably read Olivier Darrigol's text, *Should the history of relativity be revised ?* (Faut-il réviser l'histoire de la relativité ?), published in the *Lettre de l'Académie des sciences*, No. 14, Paris, 2004.

51. *In the second place I shall suppose that the forces between uncharged particles, as well as those between such particles and electrons, are influenced by a translation in quite the same way as the electric forces in an electrostatic system.*

Poincaré will also return to this question in his text *On the Dynamics of the Electron*, he writes (*op. cit.* p. 166) :

> « *Thus, Lorentz's theory would completely explain the impossibility of demonstrating absolute motion, if all forces were of electromagnetic origin. But there are forces whose origin cannot be attributed to electromagnetism, such as gravitation. Indeed, it can happen that two systems of bodies produce equivalent electromagnetic fields, that is, exerting the same action on electrified bodies and on currents, and yet these two systems do not exert the same gravitational action on Newtonian masses. The gravitational field is therefore distinct from the electromagnetic field. Lorentz was therefore obliged to complete his hypothesis by assuming that forces of all origins, and in particular gravitation, are affected by a translation (or, if one prefers, by the Lorentz transformation) in the same way as electromagnetic forces.* »

Thus, Lorentz and Poincaré had forged the mathematical key to the new mechanics : the Lorentz transformations. But they did so in service of a profoundly conservative, even **neo-Aristotelian**, principle. Their goal was to preserve the existence of an absolute reference frame—the immobile aether. Their transformations were seen as describing the curious physical effects (like contraction) that this aether imposed on moving bodies. It would fall to Einstein to propose a more radical solution : that the key fit a different lock entirely, one where the aether did not exist at all.

There is therefore truly an important question of principle here : Lorentz transformations must apply in the same way to all forces to which matter is subject, regardless of their nature, if only for the sake of consistency with observation. This is also Einstein's conclusion, who chooses this premise as the basis of his theory of relativity [52] (*op. cit.* p. 891) :

> « *In all coordinate systems for which the equations of mechanics hold, the equivalent electrodynamic and optical laws are also valid* […] *In the following,* **we raise this conjecture to the rank of a postulate** *(which we will call henceforth "principle of relativity").* »

This is the originality of Einstein's relativity, which is technically a rupture with Galilean mechanics. In his pedagogical work *Relativity : The Special and General Theory* [53], — § The Heuristic Value of the Theory of Relativity — he states his 1905 principle more concisely [54] :

> « **The general laws of nature are co-variant with respect to Lorentz transformations**. »

As we have already said, the consequences of his principle of relativity on the nature of time do not escape him. This is what he specifies, again in *Relativity* — § Minkowski's Four-Dimensional Space. Indeed, we read there :

> « *If we have not been accustomed to look at the world in this sense, as a four-dimensional continuum, it is because* **in physics**, *before the advent of the theory of relativity,* **time played a different and more independent role, compared to space**. *[...] It is for this reason that we have been accustomed to treat time*

---

52. *Wir wollen diese Vermutung (deren Inhalt im folgenden „Prinzip der Relativität" genannt werden wird) zur Voraussetzung erheben.*

53. Albert Einstein. *Relativity : The Special and General Theory* (La Relativité). Payot, collection *Petite bibliothèque*, Paris, 1990 — Original version : *Über die spezielle und die allgemeine Relativitätstheorie*, Verlag von Johann Ambrosius Barth, Leipzig 1916.

54. This is actually referring to the Poincaré group here.

> *as an independent continuum. And **it is a fact, that for classical mechanics, time is absolute, i.e., independent of position and state of motion**.* »

As we will see in a technically more precise way, the introduction of these transformations is incompatible with the notion of a universal simultaneity. With Einstein, the last pillar of the Aristotelian world—absolute Time—is demolished. The intuitive notions of a universal 'now' and the simple synchronization of clocks are revealed to be untenable, doing justice to Maimonides' ancient intuition. We are left with a unified, four-dimensional spacetime, but its geometry is no longer the simple product of Euclidean Space and Time. The task is now to construct the Inertia Group that preserves this new, strange reality.

## 10. The Poincaré Group

The passage from Galilean to Einsteinian mechanics is more dramatic than the first rupture. In the shift from Aristotle, Space was the scapegoat, but Time, the absolute, universal clock, survived. Now, in the final rupture, nothing of the old world remains. Both Aristotelian Space and Aristotelian Time are lost, dissolved into a unified four-dimensional continuum whose only remaining structure is an affine geometry endowed with a mysterious, non-positive quadratic form.

This is indeed the immediate result of the fundamental nature of Lorentz transformations, whose nature Poincaré specifies in his article as being [55]

> « …/…*a linear transformation that does not alter the quadratic form*
>
> $$x^2 + y^2 + z^2 - t^2.$$ »

---

55. H. Poincaré, *Sur la dynamique de l'électron*, Rendiconti del circolo matematico di Palermo, vol. 21, p. 146, 1906. Units are chosen such that the speed of light $c$ is 1.

It is with this proposition that the new mechanics is expressed for the first time, in such a radical and formal way. From Aristotle's Space and Time, an affine structure persists on a mixed 4-dimensional *spacetime*. The Euclidean structures of the sheets of simultaneous events and the measure of Time vanish, and are replaced by Poincaré's quadratic form. It is this space, endowed with this quadratic form, that is commonly called **Minkowski Spacetime**, we denote it ET.

From the affine structure, the new mechanics retains the most general possible notion of unforced motion : the uniform rectilinear motions, i.e., the affine lines drawn in spacetime. Thus, the new axiom of inertia is :

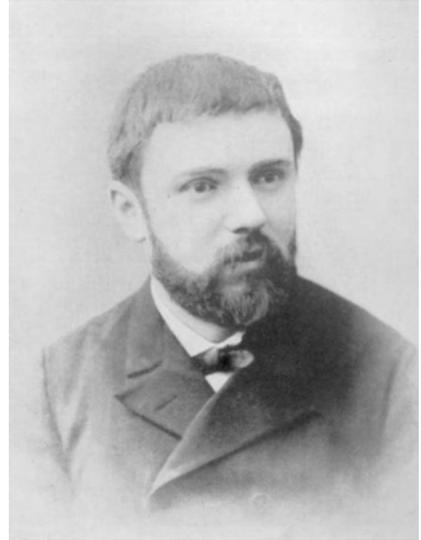

***Henri Poincaré***
*Nancy 1854, Paris 1912*

**Definition.** *All affine lines of spacetime are the* inertial motions *of Einsteinian mechanics.*

From Poincaré's quadratic form, we infer the main characteristic of the Inertia Group of Einsteinian relativity [56]. An *inertial automorphism, i.e.*, an element of the Inertia Group of Einsteinian mechanics, is an automorphism of spacetime ET that satisfies the following conditions :

- ✔ **Preserves inertial motions.** It transforms any affine line in spacetime into another affine line.
- ✔ **Preserves Poincaré's quadratic structure**. The product of two vectors by the quadratic form is unchanged under the action of this element.

---

56. To be precise, one should speak of the Einstein-Poincaré theory of relativity, or even Einstein-Lorentz-Poincaré. We choose to keep it short.

By choosing the canonical model of Minkowski spacetime, sometimes denoted $\mathbf{R}^{3,1}$, this group is called the *Poincaré Group*.

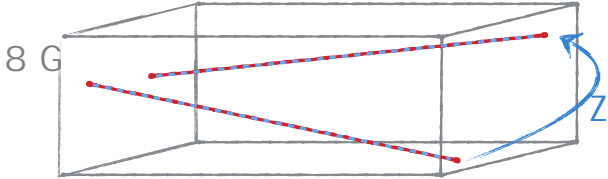

Minkowski Spacetime and the Poincaré Group.

It is described by the following proposition.

**Proposition** *The inertial automorphisms of Einsteinian relativity are the affine maps*

$$g : X \mapsto LX + C,$$

*where* $X \in ET$, $C \in \mathbf{R}^4$, *and* L *is a linear automorphism, a Lorentz transformation. Setting* $X = (r, t) \in \mathbf{R}^3 \times \mathbf{R}$, *any element of the Lorentz group decomposes uniquely as the product of two matrices*

$$L = \begin{pmatrix} \mathbf{1}_3 & \boldsymbol{\beta} \\ \bar{\boldsymbol{\beta}} & 1 \end{pmatrix} \begin{pmatrix} B^{-1} & 0 \\ 0 & a \end{pmatrix} \quad with \quad \bar{B}B + \boldsymbol{\beta}\bar{\boldsymbol{\beta}} = \mathbf{1}_3 \quad and \quad a = \pm\frac{1}{\sqrt{1-\beta^2}}.$$

B *is a* $3 \times 3$ *matrix,* $\boldsymbol{\beta}$ *is a vector in* $\mathbf{R}^3$, *the bar above* B *or* $\boldsymbol{\beta}$ *denotes the transpose, and* $\beta = \|\boldsymbol{\beta}\|$.

*The projection* $L \mapsto \boldsymbol{\beta}$ *is a trivial fibration of the* Lorentz Group $\mathfrak{L}$ *over* $\mathbf{R}^3$ *with fiber the subgroup* $O(3) \times \{\pm 1\}$. *The Poincaré group is therefore a six-dimensional group with four connected components.*

NOTE. *As for the groups of Aristotle and Galileo, it is customary to consider only the identity component of the Poincaré group.*

*Proof.* Let L be a $4 \times 4$ matrix, and G be the Gram matrix of the Poincaré quadratic form. L represents a Lorentz transformation if $\bar{L}GL = G$. Decompose $\mathbf{R}^4$ into $\mathbf{R}^3 \times \mathbf{R}$ and write, in matrix notation, this gives :

$$L = \begin{pmatrix} A & b \\ \bar{c} & a \end{pmatrix}, \quad \text{and} \quad \begin{pmatrix} \bar{A} & c \\ \bar{b} & a \end{pmatrix} \begin{pmatrix} \mathbf{1}_3 & 0 \\ 0 & -1 \end{pmatrix} \begin{pmatrix} A & b \\ \bar{c} & a \end{pmatrix} = \begin{pmatrix} \mathbf{1}_3 & 0 \\ 0 & -1 \end{pmatrix},$$

where A is a $3 \times 3$ matrix, $b, c \in \mathbf{R}^3$, and $a \in \mathbf{R}$. Thus we have the equations

$$\left. \begin{array}{ll} \bar{A}A - c\,\bar{c} = \mathbf{1}_3, & \bar{A}b - a\,c = 0 \\ \bar{b}A - a\,\bar{c} = 0, & \bar{b}\,b - a^2 = -1. \end{array} \right\}$$

Since $a \neq 0$ is evident, we have

$$c = \frac{\bar{A}b}{a}, \quad a = \pm\sqrt{1 + \bar{b}\,b} \quad \text{and} \quad \bar{A}A - \bar{A}\frac{b}{a}\frac{\bar{b}}{a}A = \mathbf{1}_3.$$

Setting then

$$\boldsymbol{\beta} = \frac{b}{a}, \quad \text{we have} \quad L = \begin{pmatrix} A & a\boldsymbol{\beta} \\ \bar{\boldsymbol{\beta}}A & a \end{pmatrix} = \begin{pmatrix} \mathbf{1}_3 & \boldsymbol{\beta} \\ \bar{\boldsymbol{\beta}} & 1 \end{pmatrix}\begin{pmatrix} A & 0 \\ 0 & a \end{pmatrix}.$$

And the condition on A and $b$ becomes $\bar{B}B + \boldsymbol{\beta}\bar{\boldsymbol{\beta}} = \mathbf{1}_3$, with $B = A^{-1}$, this is what is announced. Next, L is uniquely defined by B and $\boldsymbol{\beta}$. For $\boldsymbol{\beta} = 0$, we have $\bar{B}B = \mathbf{1}_3$ and $a^2 = 1$, *i.e.*, the subgroup $O(3) \times \{\pm 1\}$. The quotient $\mathcal{L}/O(3) \times \{\pm 1\}$ is therefore the set of eligible $\boldsymbol{\beta}$. But for any $\boldsymbol{\beta} \in \mathbf{R}^3$ with $\|\boldsymbol{\beta}\| < 1$, one can always find a B that satisfies the required condition. It suffices to take

$$B = \begin{pmatrix} \sqrt{1 - \beta^2} & 0 \\ 0 & \mathbf{1}_2 \end{pmatrix}, \quad \text{in a base such that} \quad \boldsymbol{\beta} = \begin{pmatrix} \beta \\ 0 \end{pmatrix}.$$

Thus any vector $\boldsymbol{\beta} \in \mathbf{R}^3$ with $\|\boldsymbol{\beta}\| < 1$ is eligible and the quotient is diffeomorphic to the open ball in $\mathbf{R}^3$. The Lorentz group is therefore indeed 6-dimensional, diffeomorphic (but not isomorphic) to the product $\mathbf{R}^3 \times O(3) \times \{\pm 1\}$. And since $O(3)$ has two connected components, $\mathcal{L}$ has four. □

More than ever with the

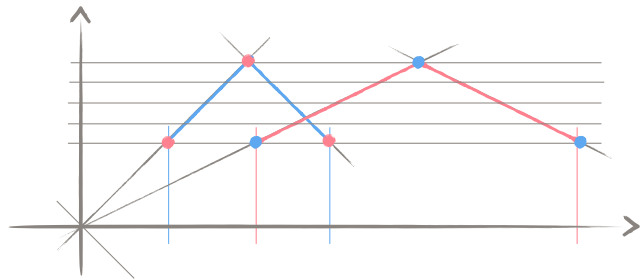

**Measuring Proper Time**

Poincaré's quadratic form is the instrument for measuring time, as the accident of the motion (Maimonide). It is called the *proper time* of the motion. For a light ray, the *proper duration* is zero over any part of the motion. Measured in any inertial embedded system, the speed of light is equal to $c$.

advent of Einsteinian relativity, mechanics is a geometry in the sense of Klein : the geometry of the Poincaré group. The intuitive (and ultimately naive) concepts of space and time, so natural to us, are dissolved in **a mysterious non-positive quadratic form whose nature** **remains improbable**. It is then the homogeneous spaces of the group of its automorphisms — the Poincaré group — that reveal the mechanics hidden behind the geometry. For example, inertial motions are partitioned into orbits of the inertia group, *space-like* motions for which the Poincaré form is positive, *time-like* motions for which the form is negative, and *light-like* motions on which the form vanishes. Here, in this particular framework of special relativity, mechanics can be nothing other than a geometry, Poincaré geometry.

## 11. Time and Simultaneity

As we have seen, the notion of Aristotelian Time and the notion of *simultaneity* are intimately linked, if not identical. We will interpret them precisely based on the preceding constructions and draw certain conclusions.

We continue to denote by ET the spacetime continuum — whether it refers to Aristotle's Space × Time, or Galilean spacetime, or Minkowski space. It is considered solely for its differential and affine structure, and is equivalent to $\mathbf{R}^4$. The Inertia Group G of the chosen mechanics, be it Aristotle's group, Galileo's group, or the Poincaré group, acts naturally on this continuum. The *simultaneity of events* can then be interpreted as a (infinitely) differentiable fibration $\tau : \mathrm{ET} \to \mathrm{T}$, where T is Aristotle's time, equivalent to the real line $\mathbf{R}$ equipped with its standard Euclidean structure, oriented from past to future by a measure of time (i.e., any positive vector). Simultaneous events, correspond to *instants* — the elements of T — and are, by definition, their preimages. In other words, two events $x$ and $x'$ are simultaneous if they belong to the same *fiber* of $\tau$, *i.e.*, if $\tau(x) = \tau(x')$. The coherence of simultaneity, with the chosen mechanics, corresponds to the equivariance of this *temporalization*, by the group of mechanics. *i.e.*, there exists a map $h$ from group G to the group of Euclidean transformations of the temporal line, such that, for any $x \in \mathrm{ET}$ and for any $g \in G$, we have :

$$\tau(g(x)) = h(g)(\tau(x)) \quad i.e. \quad \tau(g(x)) = \varepsilon(g)\tau(x) + \varphi(g),$$

where $\varepsilon(g) = \pm 1$ and $\varphi : G \to \mathbf{R}$. This simply means that *simultaneous events are transformed into simultaneous events by the group of mechanics*, and that this transformation depends only on the group [57]. Considering only orthochronous transformations : those that do not reverse the course of time, *i.e.*, the direction (orientation) of motions, we have $\varepsilon(g) = +1$ for all $g$. Then, since group G acts on ET, *i.e.*, $g(g'(x)) = (gg')(x)$, we deduce that the map $\varphi$, which is the action of the group of mechanics on T, is a (smooth) homomorphism from G to the additive group of real numbers : $\varphi(gg') = \varphi(g) + \varphi(g')$. And since Aristotelian time flows indefinitely from past to future, this homomorphism is necessarily surjective [58].

This formalization is the key to the entire argument. It transforms a vague philosophical debate about the 'absoluteness' of time into a precise, mathematical question : does such a non-trivial map from the Inertia Group to the real line exist ? The answer, which depends entirely on the structure of the group, will be definitive.

**Statement** *An Aristotelian time is necessarily associated, by the foregoing, with a surjective smooth homomorphism from the Inertia Group of mechanics, to the additive group of real numbers.*

As shown by the following theorem [59], we immediately deduce the mathematical translation of the absolute incompatibility between the existence of an Aristotelian Time and Einsteinian relativity. Stated differently, it is the structural impossibility of defining a universal

---

57. This is in particular what Einstein means when he says (see above) « *absolute time [is] independent of position and state of motion* ».

58. We have implicitly assumed oriented spaces and group actions preserving orientation.

59. This theorem relies on only one slightly sophisticated result : the Lorentz group is simple. The reader can find a proof in :

Eugene Wigner. *On Unitary Representations Of The Inhomogeneous Lorentz Group.* Annals of Mathematics, vol. 40 (1), pp. 149–204, 1939.

simultaneity of events compatible with the principle of relativity—or, as physicists would say, the impossibility of synchronizing clocks throughout spacetime in a manner that is consistent with Poincaré invariance.

**THEOREM.** *There exists no smooth homomorphism, from the Poincaré group onto the group of translations of the real line.*

*Proof.* Let $\varphi : G \to \mathbf{R}$ be a homomorphism from the Poincaré group G to the real line equipped with addition. Consider the groups $\mathfrak{L}$ and $(\mathbf{R}^4, +)$, injected into G as subgroups, in the following way :

$$\begin{pmatrix} L & 0 \\ 0 & 1 \end{pmatrix} \text{ with } L \in \mathfrak{L} \quad \text{and} \quad \begin{pmatrix} \mathbf{1} & C \\ 0 & 1 \end{pmatrix} \text{ with } C \in \mathbf{R}^4,$$

and denote $\phi$ and $f$, the restrictions of $\varphi$ to $\mathfrak{L}$ and $(\mathbf{R}^4, +)$,

$$\phi(L) = \varphi \begin{pmatrix} L & 0 \\ 0 & 1 \end{pmatrix} \quad \text{and} \quad f(C) = \varphi \begin{pmatrix} \mathbf{1} & C \\ 0 & 1 \end{pmatrix}.$$

The restrictions $\phi$ and $f$ are therefore (smooth) homomorphisms from $\mathfrak{L}$ and $(\mathbf{R}^4, +)$ to $(\mathbf{R}, +)$. Decomposing the elements of G in the following way :

$$\begin{pmatrix} L & C \\ 0 & 1 \end{pmatrix} = \begin{pmatrix} \mathbf{1} & C \\ 0 & 1 \end{pmatrix} \begin{pmatrix} L & 0 \\ 0 & 1 \end{pmatrix}$$

we obtain

$$\varphi \begin{pmatrix} L & C \\ 0 & 1 \end{pmatrix} = f(C) + \phi(L). \tag{$\clubsuit$}$$

But $f(C) + \phi(L) = \phi(L) + f(C)$, *i.e.* :

$$\phi(L) + f(C) = \varphi \begin{pmatrix} L & 0 \\ 0 & 1 \end{pmatrix} + \varphi \begin{pmatrix} \mathbf{1} & C \\ 0 & 1 \end{pmatrix} = \varphi \begin{pmatrix} L & LC \\ 0 & 1 \end{pmatrix} = f(LC) + \phi(L).$$

So $f(LC) = f(C)$, *i.e.* : $f$ is a homomorphism from $(\mathbf{R}^4, +)$ to $\mathbf{R}$, invariant under the Lorentz group. Conversely, it is immediate to verify that any pair $(f, \phi)$ of such homomorphisms defines, by $(\clubsuit)$, a homomorphism $\varphi : G \to \mathbf{R}$, and this decomposition is unique.

Next, since $f$ is invariant under the Lorentz group, $f(C) = F(|C|)$, for some F and where $|\cdot|$ denotes the Poincaré quadratic form. But $f$ is a homomorphism, so $F(|C + C'|) = F(|C|) + F(|C'|)$ (an identity which cannot hold for a non-zero function, since the Minkowski quadratic form is not linear ; for example, for two non-collinear space-like vectors, the 'triangle inequality' is reversed). This immediately implies

$$f = 0 \quad \text{and thus} \quad \varphi\begin{pmatrix} L & C \\ 0 & 1 \end{pmatrix} = \phi(L).$$

Let then $K \subset \mathfrak{L}$, be the kernel of the homomorphism $\phi$, and $H = \phi(\mathfrak{L}) \subset \mathbf{R}$ its image. The group $K$ is a normal subgroup of the Lorentz group, which is a simple group [60], thus $K = \{\mathbf{1}\}$ or $\mathfrak{L}$. The homomorphism $\phi$ being differentiable, considering its tangent linear map at the identity, we have : $\dim(K) + \dim(H) = \dim(\mathfrak{L})$. Since $\dim(H) = 0$ or $1$, and $\dim(K) = 0$ or $\dim(\mathfrak{L}) = 6$, the only possibilities remaining are $\dim(K) = 6$ and $\dim(H) = 0$, *i.e.*, $\phi = 0$. Thus $\varphi = 0$, and there is no surjective homomorphism from the Poincaré group onto the real line. $\qquad\qquad\square$

**Remark 1.** The notion of *instant*, as *simultaneity* of events, expresses the difference between time as a chronological parameter, an accident of motion — which persists in relativity — and Aristotle's Time, absolute, which intends to be a unique and coherent reference for all motions. The former is at the source, in the description of motion, it is an *accident* of motion, it accompanies it, it is the parameter of its event chronology. The latter is (or is not) at the goal of the fibration by simultaneity, it defines the absolute instant, true for all motions of all things, it is Aristotle's Time.

**Remark 2.** The preceding theorem shows unequivocally the futility of neo-Aristotelian attempts to restore a universal time, in mechanics or cosmology, which one can sometimes encounter in the literature.

Once again, it is observed, that any statement concerning physics, once clearly expressed within the formal, rigorous, framework of mathematics — here in particular geometry — finds its ultimate, indisputable expression. It provides both an unambiguous answer to the question concerned, but also clearly delineates the domain of validity of the statement and its correspondence with the real world."'

---

60. See Wigner, *op. cit.* section D, p.167.

## 12. The Three Historical Ruptures in Mechanics

The journey through the history of motion according to place has revealed three capital epistemological ruptures. Each represents the collapse of a geometry and the construction of a new one, built upon a more subtle axiom of inertia.

**The first rupture** is a hinge, the *initialization*, it is the establishment of the categories of Time and Space by Aristotle, primary categories through which motion is naturally expressed in a tangible way. The universal, eternal, and absolute character of the categories of Space and Time is intimately dependent on the immobility of the Earth at the center of the universe, as well as the perfection (homogeneity) of the spheres' motion and their eternity [61].

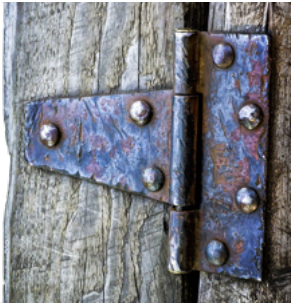

*A **Hinge**, both rupture and link.*

**The second rupture**, we owe it to Giordano Bruno and Galileo Galilei; it is the abandonment of Space, which is not what it claims to be. Rest being indistinguishable from uniform rectilinear motion, Space as the place of rest for things is no longer founded. We are here confronted with a complex phenomenon, because if Space is obsolete, it is not immediately replaced; no one realizes or fully understands how to handle this situation, and Aristotle's Space continues to linger in the minds of physicists as if nothing had happened.

**The third rupture** is more radical, because by putting an end to Aristotelian Time, it also destroys everything that remained of our intuition about the decomposition of motion into Time and Space. Of the original spacetime, only a global affine structure equipped with a quadratic form remains. This compels us to introduce a new paradigm, and recognize this new mechanics as the geometry associated with the group of automorphisms of this quadratic form. Einsteinian

---

61. See for example *The Immobile Earth* by Jean-Jacques Szczeciniarz, PUF, 2003.

relativistic mechanics becomes Poincaré geometry, in the sense of Felix Klein.

**The reconstruction** All this necessarily leads us to rethink the nature of the previous iterations : Aristotelian and Galilean mechanics ; to reinterpret them, under penalty of inconsistency, as geometries, *i.e.*, to exhibit the group from which they emanate. This is what we have done :

- ⇢ We have assumed the existence of a 4-dimensional affine spacetime, common to all three cases.
- ⇢ For each case, we have singled out a family of motions that we have called *inertial motions*, in this case, a family of lines in spacetime.
- ⇢ We have specified for each case a structure, foliations or quadratic form, that characterizes it.

We then obtained the *Inertia Groups* of these mechanics as the groups of automorphisms of these structures, preserving the family of inertial motions. This gave us :

- ✔ The Group of Aristotle.
- ✔ The Galilean Group.
- ✔ The Poincaré Group.

In return, each of these groups restored the mechanics associated with it as a geometry in the sense of Klein. The structures from which these groups arise can then fade behind the group itself, since they are fully embedded within it. There is in the process we have just described an enrichment that could be called dialectical : the group sublimes the structure from which it originates. We will see some particular consequences later, in the context of symplectic mechanics.

This geometric approach, centered on the Inertia Group, is more than a matter of mathematical taste ; it represents a fundamental choice in the philosophy of physics. Dynamical laws, such as Newton's second law, are powerful computational tools. Yet the very essence of dynamics lies in a simple act of comparison : a given motion is understood by measuring its deviation from a privileged family of reference motions. And this is what dynamics will deliver. For Aristotle, this family was that of the resting motions; one measures movement against immobility. More generally, the set of *inertial motions* constitutes the fundamental basis against which all other motions are analyzed.

The geometry of a mechanics is then, by definition, the group of spacetime automorphisms that preserves this basis, exchanging one inertial motion for another. This group is not merely a property of the mechanics; it *is* the mechanics, defining its character and its limits. Rephrasing Newton : just as the compass, a mechanical device, gives birth to the circle, so too does a specific principle of inertia give birth to a geometry. In this view, Euclidean geometry is nothing other than the formal expression of Aristotelian mechanics.

**Primary and secondary ruptures** Meanwhile, we have seen, thanks to this construction, how an often vague epistemological discourse — on the evolution of concepts in mechanics since Aristotle — crystallizes into a concise, precise mathematical object, suffering no ambiguity, and itself perfectly discriminating as to the epistemological rupture it represents. We can now propose a rigorous hierarchy for these historical shifts : a **primary rupture** is one that fundamentally alters the Inertia Group of mechanics, while a **secondary rupture** introduces new dynamical laws *within* a pre-existing, stable geometric framework. One can think in particular of the Newtonian rupture. It is a hinge between the sublunar world and the supralunar world. Newton teaches us, contrary to what was believed, that it is the same physics that governs these two worlds. This is obviously a profound rupture in cosmogony and mechanics, but it is nevertheless a **secondary** or **subordinate rupture**. It does not question the mechanics' group, whose domain extends but remains that of Galilean geometry. We thus observe a hierarchy solidly defined by geometry, in the nature of the epistemological ruptures that weave the history of mechanics.

**Historically** Finally, let us note that these three major eras that we have evoked, where mechanics was seen to evolve from Aristotelian concepts of Space and Time to their abandonment, are spread over periods of rapid decrease, marking a spectacular acceleration in the history of ideas in mechanics, as shown in the following diagram.

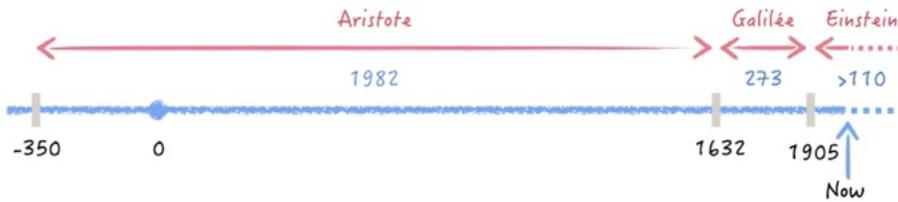

The acceleration of ideas in Mechanics.

The next volume of this study will focus precisely on *The Analytical Revolution* in mechanics, which was initiated by Newton with his memoir on the motion of planets, which later became the *Principia*. He there developed the notion of force, or more precisely, he described its measure, by the variation from uniform rectilinear motion—the first great act in the **art of comparing real motions to ideal ones**, which is the very definition of dynamics. We will then see how Lagrange renewed the study of natural dynamical systems by revealing a sophisticated structure on the spaces of motions, which we now call a symplectic structure. In a certain sense, which we will develop, this structure has become the operational response to the replacement of the obsolete categories of Space and Time.

# 13. Bibliographical References

# 14. Suggested Readings

* Plato, *Timaeus*, transl. Émile Chambry.

* Aristotle, *Physics*, Ed. J. Vrin.

* Maimonides, *The Guide for the Perplexed*, Ed. Maisonneuve & Larose.

* Giordano Bruno, *The Ash Wednesday Supper*, Ed. L'éclat.

* Galileo Galilei, *Dialogue Concerning the Two Chief World Systems*, Ed. Points.

* Joseph-Louis Lagrange, *Analytical Mechanics*, Ed. Blanchard.

* Albert Einstein, *Relativity : The Special and General Theory*, Ed. Payot.

* Felix Klein, *The Erlangen Program*, Ed. Gauthier-Villars.

* Jean-Marie Souriau, *Structure of Dynamical Systems* (Structure des Systèmes Dynamiques), Éd. Dunod.

About the author

piz@math.huji.ac.il          http://math.huji.ac.il/~piz

Patrick Iglesias-Zemmour is currently at the Hebrew University of Jerusalem (Israel), after having been a researcher at the CNRS (France).

# THE GEOMETRY OF MOTION
## PART I : MECHANICS AS GEOMETRIES

---

### PATRICK IGLESIAS-ZEMMOUR

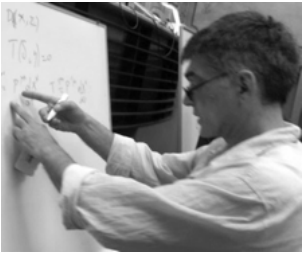

At the beginning of Physics lies motion, for as Aristotle writes, "Nature is a principle of motion and change." Among these changes, our focus is on "motion according to place." Aristotle was the first to provide its formal framework by introducing the intuitive categories of Time and Place—what we now call Space.

Yet, the discoveries that followed, from the ancients to the moderns, from Aristotle to Einstein by way of Galileo, forced us to abandon these categories. These abandonments came at a cost, for they contradicted our deepest intuitions about the world. We speak of the "loss of Space" or the "end of absolute Time," but what do these phrases truly mean? Can we precisely locate and capture the epistemological ruptures that led us from the physics of Aristotelian equilibrium to the dynamics of Einstein's relativity?

This book answers with a definitive yes.

This first volume in an exploration of The Geometry of Motion identifies the three great ruptures in the history of mechanics. It demonstrates how to associate each one —Aristotelian, Galilean, and Einsteinian— with a unique Inertia Group. This group of transformations acts on the common substrate of spacetime, and it is this group that captures the absolute essence of the mechanics. By translating physical principles into the language of group theory, this book reveals each mechanics for what it truly is : a distinct and coherent Geometry in the sense of Felix Klein.